\documentclass[reqno]{amsart} 
\usepackage{epsfig,amscd,amsmath,amssymb,amsthm}
\usepackage[all]{xy}
\newtheorem{thm}{Theorem}[section]

\newtheorem{lem}[thm]{Lemma}

\newtheorem{prop}[thm]{Proposition}

\theoremstyle{definition}
\newtheorem{defn}[thm]{Definition}
\newtheorem{exmp}[thm]{Example}
\newtheorem{rem}[thm]{Remark}

\newtheorem{alg}{Algorithm}
\long\def\@makealgocaption#1#2{\vskip 2ex \small
  \hbox to \hsize{\parbox[t]{\hsize}{{\bfseries #1.} #2}}}
\newcounter{algorithm}
\def\thealgorithm{\@arabic\c@algorithm}
\newtheorem{examples}[thm]{Examples}

\newcommand{\A}{{\mathbb A}}

\newcommand{\C}{{\mathbb C}}

\renewcommand{\P}{{\mathbb P}}

\newcommand{\K}{K}   %----> ground field
\newcommand{\kd}{{\mathcal D}}

\newcommand{\kf}{{\mathcal F}}
\newcommand{\kg}{{\mathcal G}}
\newcommand{\kh}{{\mathcal H}}

\newcommand{\ko}{{\mathcal O}}

\newcommand{\ks}{{\mathcal S}}

\newcommand{\fm}{\mathfrak{m}}

\DeclareMathAlphabet{\mathsc}{U}{rsfs}{m}{n}
\newcommand{\sA}{\mathsc{A}}
\newcommand{\sC}{\mathsc{C}}

\newcommand{\sM}{\mathsc{M}}

\DeclareMathOperator{\charac}{char}
\DeclareMathOperator{\Sing}{Sing}
\DeclareMathOperator{\Spec}{Spec}

\newcommand{\ba}{\boldsymbol{a}}

\newcommand{\bb}{\boldsymbol{b}}

\newcommand{\beps}{\boldsymbol{\varepsilon}}

\newcommand{\bo}{\boldsymbol{0}}

\newcommand{\bp}{\boldsymbol{p}}

\newcommand{\bs}{\boldsymbol{s}}

\newcommand{\Ies}{I^{\text{\it es}}}
\newcommand{\IES}{I^{\text{\it ES}}}
\newcommand{\Iesf}{I^{\text{\it es}}_{\text{\it fix}}}
\newcommand{\Iesphi}{I^{\text{\it es}}_{\varphi}}
\newcommand{\Ovsigma}{\overline{\sigma}}
\newcommand{\Ovkc}{\overline{\sC}}

\newcommand{\OvO}{\overline{0}}

\newcommand{\Keps}{K[\varepsilon]}

\newcommand{\lra}{\longrightarrow}

\newcommand{\osigma}{\overline{\sigma}}

\newcommand{\fa}{{\mathfrak a}}
\newcommand{\fb}{{\mathfrak b}}

\DeclareMathOperator{\codim}{codim}

\DeclareMathOperator{\Def}{\kd\!\!\:\it{ef}}

\DeclareMathOperator{\Res}{Res}

\newcommand{\uDef}{\underline{\Def}\!\,}

\newcommand{\Defes}{\Def^{\text{\it es}}}
\newcommand{\DefES}{\Def^{\text{\it ES}}}

\newcommand{\uDefes}{\uDef^{\text{\it es}}}
\newcommand{\uDefES}{\uDef^{\text{\it ES}}}

\newcommand{\Defsec}{\Def^{\text{\it sec}}}
\newcommand{\Defsecii}{\Def^{\text{\it sec}}}
\newcommand{\uDefsec}{\uDef^{\text{\it sec}}}
\newcommand{\uDefsecii}{\uDef^{\text{\it sec}}}

\DeclareMathOperator{\Der}{Der}

\DeclareMathOperator{\ES}{\it ES}

\DeclareMathOperator{\id}{id}

\DeclareMathOperator{\modulo}{mod}
\renewcommand{\mod}{\text{mod } }

\DeclareMathOperator{\ord}{ord}

\DeclareMathOperator{\pr}{pr}
\DeclareMathOperator{\Quot}{Quot}

\newcommand{\Rbar}{\overline{R}}

\DeclareMathOperator{\Sets}{\ks\it{ets}}

\newcommand{\Ties}{T^{1,\text{\it es}}}
\newcommand{\TiES}{T^{1,\text{\it ES}}}
\newcommand{\Tisec}{T^{1,\text{\it sec}}}
\newcommand{\Tisecii}{T^{1,\text{\it sec}}}
\newcommand{\koC}{\ko_{C,0}}

\newenvironment{entry}
  {\begin{list}{}%
     {%
      \setlength{\labelwidth}{45pt}%
      \itemsep1pt%
      \setlength{\leftmargin}{53pt}%
      }%
  }%
{\end{list}}

\begin{document}

\title{Equisingular Calculations for Plane Curve Singularities}

\author{Antonio Campillo}
\address{Departamento de Algebra, Geometria y Topologia\\
Universidad de Valladolid\\Facultad de Ciencias\\
E -- 47005   Valladolid}
\author{Gert-Martin Greuel}
\address{
Fachbereich Mathematik\\
TU Kaiserslautern\\
Erwin-Schr\"odinger-Stra{\ss}e\\
D -- 67663 Kaiserslautern}
\author{Christoph Lossen}
\address{
Fachbereich Mathematik\\
TU Kaiserslautern\\
Erwin-Schr\"odinger-Stra{\ss}e\\
D -- 67663 Kaiserslautern}

\begin{abstract}
We present an algorithm which, given a deformation with section of a reduced
plane curve singularity, computes equations for the equisingularity stratum
(that is, the $\mu$-constant stratum in characteristic $0$) in the
parameter space of the deformation. The algorithm works for any, not
necessarily reduced, parameter space and for algebroid curve
singularities $C$ defined over an algebraically closed field of
characteristic $0$ (or of characteristic \mbox{$p>\ord(C)$}). It
provides at the same time an algorithm for computing the
equisingularity ideal of J.\ Wahl. The algorithms have been
implemented in the computer algebra system {\sc Singular}. We show
them at work by considering two non-trivial examples.
As the article is also meant for non-specialists in singularity
theory, we include a short survey on new methods and results about
equisingularity in characteristic $0$.
\end{abstract}

\maketitle

\centerline{\small \sl Dedicated to the memory of Sevin Recillas}
\bigskip

%\begin{keyword}
%  Equisingular deformation, plane curve singularities, Milnor number,
%  Ham\-bur\-ger-Noe\-ther expansion, deformation of parametrization
%\end{keyword}

\section{Introduction}
\noindent
Equisingular families of plane curve singularities, starting from Zariski's
pioneering 'Studies in Equisingularity I--III' \cite{Zar}, have been of
constant interest ever since. Zariski intended to develop this concept aiming
at a resolution of singularities where 'equisingular' singularities should
resolve simultaneously or are, in some sense, natural centres for blowing
up. This approach was completely successful only in the case of families of
plane curves\footnote{Zariski originally considered equisingularity of a
  (germ of a) hypersurface $X$ along a subspace \mbox{$Y\subset X$} and a projection
  of \mbox{$\pi:X\rightarrow T$} such that $Y$ is the image of a section of
  $\pi$. If $Y$ has codimension 1 then the fibres of \mbox{$X\to T$} are plane curve
  singularities. Zariski then considered the discriminant of the projection
  which is a hypersurface in $T$ (at least if $T$ is smooth) and thus
  equisingularity of $X$ along $Y$ can be defined by induction on the
  codimension of $Y$ in $X$.} where Zariski introduced several quite
different, but equivalent, notions of equisingularity.

One of these notions was used by J.\ Wahl in his thesis to extend
the concept of equisingularity to families over possibly non-reduced base
spaces (see \cite{Wah}). This enabled him to apply Schlessinger's
theory of deformations over 
Artinian rings and to define the equisingularity ideal wich describes the
tangent space to the functor of equisingular deformations. Moreover, Wahl
proved that the base space of the semiuniversal equisingular deformation of a
reduced plane curve singularity is smooth. Wahl's proof of this theorem, which
is an important result in singularity theory, is quite complicated and uses
several intermediate deformation functors, in particular deformations of the
exceptional divisor of the embedded resolution of the singularity. Hence, he has
to pass to deformations of global objects (the exceptional divisor) and not
just of singularities.

The definition of equisingularity is algebraic and uses the resolution of
singularities. But there is also a purely topological definition: two
reduced plane curve singularities 
$(C_1,0)$ and $(C_2, 0)$ in $(\C^2, 0)$ are equisingular iff they
have the same {\em embedded topological type\/}, that is, there exist (arbitrary small)
balls \mbox{$B_1, B_2\subset \C^2$} centred at $0$ and a homeomorphism of the triple
\mbox{$(B_1, C_1\cap B_1, 0)$} onto \mbox{$(B_2, C_2\cap B_2, 0)$} for
representatives $C_i$ 
of $(C_i, 0)$. As \mbox{$(B_i, C_i\cap B_i, 0)$} is homeomorphic to the cone over
\mbox{$(\partial B_i, C_i\cap \partial B_i)$}, the topological type of a reduced plane
curve
singularity $(C,0)$ is determined by the embedding of the link \mbox{$C\cap \partial
B$} in $\partial B$,  which consists of $r$ knots (circles $S^1$ embedded in
\mbox{$\partial B\approx S^3)$} where $r$ is the number of irreducible
components of $(C,0)$. 

The topological type of each knot, which is an iterated torus knot, is
determined by the pairs of ''turning numbers'' for each iterated torus which
itself are determined by and determine the sequence of Puiseux pairs of the
corresponding branch. Moreover, the linking number of two knots coincides with
the intersection number of the corresponding two branches. Hence, the
topological type of $(C,0)$ is characterized by the Puiseux pairs of each
branch and by the pairwise intersection numbers of different branches.
This shows that the system of Puiseux pairs and the intersection numbers form a
complete set of numerical invariants for the topological type or the
equisingularity type of a plane curve singularity. 

If we consider not just
individual singularities but families, then the situation is even more
satisfactory: the topological type is controlled by a single number, the
Milnor number. Indeed, we have the following result due to Zariski \cite{Zar},
L\^{e} \cite{LDT,LeR} and Teissier \cite{Tei1}. Let \mbox{$\pi:(\sC, 0)\to
(T,0)$} be a flat family of reduced plane curve singularities with section
\mbox{$\sigma:(T,0)\to (\sC,0)$}, then the following are equivalent (for
\mbox{$\sC\to T$} a small representative of $\pi$ and
\mbox{$\sC_t=\pi^{-1}(t)$} the fibre over \mbox{$t\in T$}):

\begin{enumerate}
\item[(1)] $(\sC, 0)\xrightarrow{\pi} (T,0)$ is equisingular along $\sigma$,
\item[(2)] the topological type of $(\sC_t, \sigma(t))$ is constant
  for $t\in T$,
\item[(3)] the Puiseux pairs of the branches of $(\sC_t, \sigma(t))$
  and the pairwise intersection multiplicities of the branches are
  constant for \mbox{$t\in T$},
\item[(4)] the $\delta$-invariant $\delta(\sC_t, \sigma (t))$ and the number
  of branches $r(\sC_t, \sigma(t))$ are constant for \mbox{$t\in T$},
\item[(5)] the Milnor number $\mu(\sC_t, \sigma(t))$ is constant for
  $t\in T$.\:\footnote{By a theorem of Lazzeri, if
    $\mu(\sC_t)=\sum_{x\in \Sing(\sC_t)} 
  \mu(\sC_t, x)=\mu(C, 0)$ for $t\in T$ then there
  is automatically a section $\sigma$ such that $\sC_t\smallsetminus\sigma(t)$
  is smooth and $\mu(\sC_t, \sigma(t))$ is constant.}
\end{enumerate}

\noindent
Recall that for a reduced plane curve singularity \mbox{$(C,
  0)=\{f=0\}\subset(\C^2, 0)$} defined by a (square-free) power series
\mbox{$f\in \ko_{\C^2, 0}=\C\{x,y\}$}, the invariants
  $\mu$, $r$, and $\delta$ are defined as follows:\label{page:Milnor}
\[
\begin{array}{lcl}
\mu(C,0) & = & \dim_\C \C\{x,y\}/\langle \frac{\partial f}{\partial
  x}, \frac{\partial f}{\partial y}\rangle\,,\\
r(C,0) & = & \text{number of irreducible factors of } f\,,\\
\delta (C,0) & = & \dim_\C \overline{\ko}_{C,0}/\ko_{C,0}\,.
\end{array}
\]
Here,  \mbox{$\ko_{C,0}=\ko_{\C^2,
  0}/\langle f\rangle$} and $\overline{\ko}_{C,0}$ is the
normalization of $\koC$, that is, the integral 
closure of $\koC$ in its total ring of fractions. Furthermore, for each
reduced plane curve singularity we have the relation (due to Milnor \cite{Mil})
\[
\mu=2\delta-r+1\,.
\]
This result was complemented by Teissier \cite{Tei}\footnote{The original
  proof of Teissier and Raynaud in \cite{Tei} has been clarified and extended
  to families of (projective) varieties in any dimension by Chiang-Hsieh and
  Lipman in \cite{CHL}.}, showing that for a normal base $(T,0)$, the
flat family \mbox{$\pi: (\sC, 0)\to (T,0)$} admits a
  simultaneous normalization iff \mbox{$\delta
    (\sC_t)=\sum_{x\in\Sing(\sC_t)}\delta(\sC_t, x)$} is constant. 

The equivalence of (1) and (4) above shows the following:
Let \mbox{$(\sC ,0)\to (T,0)$} be the seminuniversal deformation of
$(C,0)$ and let 
\[
\Delta^\mu=\{t\in T\ |\ \mu(\sC_t) = \mu (C,0)\}
\]
be the $\mu$-constant stratum of $(C,0)$. Then  $\Delta^\mu$ coincides (as a set)
with the equisingularity stratum of Wahl and, hence, is smooth.

Note that for higher dimensional isolated hypersurface singularities the
$\mu$-constant stratum is in general not smooth, cf.\ \cite{Lue}.

Despite the fact that the equisingularity stratum admits such a simple
description, all attempts to find a general simple proof for its smoothness
failed (except for irreducible germs, cf.\ \cite{Tei}).

One purpose of this paper is to report on a simple proof of Wahl's
theorem. The idea is to consider deformations of the parametrization
$$\varphi:(\overline{C}, \overline{0})\to (C,0)\hookrightarrow(\C^2, 0)$$ of $(C,0)$, where
\mbox{$(\overline{C}, \overline{0})\to (C,0)$} is the normalization of $(C,0)$. We 
define equisingular deformations of  $\varphi$ and
prove that they are unobstructed. This is very easy to see, as they
are (in certain coordinates) 
even linear. Then we show (by a direct argument on the tangent level)
that equisingular deformations of $\varphi$ and 
equisingular deformations of $(C,0)$ have isomorphic semiuniversal objects. 

This proof has been known by the second author since about fifteen years and was
communicated at several conferences. A preliminary preprint \cite{GRec}, joint
with Sevin Recillas, has even been cited by some authors. Later on, these results
have been extended to positive characteristic in a joint preprint
of the authors \cite{CGL} where, in addition, an algorithm  
to compute the equisingularity stratum was developed and used to prove
one of the main results.
However, meanwhile the theory of equisingularity in positive
characteristic was further developed by the authors where the
algorithm itself could be eliminated in the theoretical arguments
\cite{CGL1}. These results will be published elsewhere, but as we
think that the algorithmic part of \cite{CGL} should not be forgotten,
we present it in this paper.

We start with a survey of the new methods and
results about equisingularity in characteristic $0$ with a sketch of
the proofs (for more details, we refer to \cite{CGL1}). The main
purpose of this paper is to describe an algorithm to compute the
$\mu$-constant stratum 
$\Delta^\mu$ for an arbitrary deformation \mbox{$(\sC,0)\to (T,0)$}
with section of
a reduced plane curve singularity $(C,0)$. More precisely, this
  algorithm computes an ideal \mbox{$I\subset \ko_{T,0}$} with
  \mbox{$\Delta^\mu=V(I)$}. As a corollary, we obtain an algorithm to
  compute the equisingularity ideal of 
  Wahl. The algorithms work also in characteristic \mbox{$p>0$} if $p$
  is larger than the multiplicity of $C$ and we
  formulate them in this generality.
They have been implemented in {\sc Singular} \cite{GPS} by
A.\ Mindnich and the third author \cite{LoM}.

\section{The Fundamental Theorems}\label{sec:theorems}
\setcounter{thm}{0}
\setcounter{equation}{0}

By Wahl, the equisingularity stratum $\Delta^\mu$ in a versal family
\mbox{$(\sC,0)\to (T,0)$}
(with section $\sigma$) is smooth. The idea of our proof for this fact is
extremely simple. Consider the parametrization 
\[
\varphi_i:(\C, 0)\lra (\C^2, 0),\quad  t_i\longmapsto (x_i (t_i), y_i(t_i))
\]
of the $i$-th branch $(C_i, 0)$ of $(C,0)$. Let, for $i=1, \ldots, r$,
\begin{equation}\label{eqn:(*)}
\begin{array}{lcl}
x_i(t_i) & = & t_i^{n_i},\\
y_i(t_i) & = & t_i^{m_i}+\displaystyle\sum\limits_{j\geq 1} a_i^j t_i^{m_i+j}.
\end{array}
\end{equation}
Now, we use the above characterization (3) for equisingularity,
assuming that $\sigma$ is the trivial section. 
Fixing the Puiseux pairs of $(C_i,0)$ is equivalent to the condition that no
new characteristic term appears if we vary the $a_i^j$. For each
$i$, this is an open condition on the coefficients $a_i^j$. Moreover, it is
easily checked that fixing the intersection multiplicity
of $(C_i,0)$ and $(C_k, 0)$ defines a linear condition among the $a_i^j$ and
$a_k^j$. Thus, if we consider \eqref{eqn:(*)} as a deformation of $(C,0)$ with $a_i^j$
replaced by coordinates $A_i^j$, \mbox{$A_i^j(0)=a_i^j$}, then the equisingular
deformations form a smooth subspace in the parameter space with coordinates
$A_i^j$. This family is easily seen to be versal. By general facts
from deformation theory it follows then that each versal equisingular
deformation of the parametrization has a smooth parameter space.

This argument works only for deformations over reduced
base spaces $(T,0)$. In particular, it does not work for {\em infinitesimal
deformations\/}, that is, for deformations over
\[
T_\varepsilon:=\Spec(\C[\varepsilon]/\langle\varepsilon^2\rangle).
\]
On the other hand, in order to use the full power of deformation theory we
need infinitesimal deformations. 

We continue this section by giving the required definitions for deformations of (the
equation of) $(C,0)$ and of the parametrization of $(C,0)$ in the framework of
deformation theory over arbitrary base spaces. These definitions are quite
technical, which is, however, unavoidable.

Throughout the following, let \mbox{$(C,0) \subset (\C^2\!,0)$} be a reduced
plane curve singularity, and let \mbox{$f \in \langle x,y\rangle^2\subset \C\{x,y\}$}
be a defining power series. We call $f=0$, or just $f$ the {\em (local)
  equation\/} of
$(C,0)$. Deformations of $(C,0)$ (respectively embedded deformations of
$(C,0)$) will also be called 'deformations of the
  equation'
(in contrast to 'deformations of the parametrization', see Definition
\ref{def:def of param}).

\begin{defn} A {\em deformation (of the equation)\/} of $(C,0)$ over
a complex germ $(T,0)$ is a flat morphism \mbox{$\phi:(\sC, 0)\to (T,0)$} of complex
germs together with an isomorphism
\mbox{$i:(C,0)\xrightarrow{\cong}(\phi^{-1}(0),0)$}. It is denoted by
$(i, \phi)$. 

A morphism from $(i,\phi)$ to a deformation \mbox{$(i', \phi'):(C,0)\hookrightarrow
(\sC',0)\to (T', 0)$} consists of morphisms \mbox{$\psi:(\sC,0)\to (\sC',0)$} and
\mbox{$\chi:(T,0)\to (T',0)$} 
making the obvious diagram commutative. If, additionally,
a section \mbox{$\sigma$} of $\phi$ is given (that
is, a morphism \mbox{$\sigma: (T,0)\to (\sC, 0)$} satisfying
\mbox{$\phi\circ\sigma=\id_{(T,0)}$}), we speak about a 
{\em deformation with section\/}, denoted by $(i, \phi, \sigma)$.
\end{defn}

\noindent
A more explicit description is as follows: since each deformation of
\mbox{$(C,0) \subset (\C^2,0)$} can be embedded, there is an isomorphism
\mbox{$(\sC, 0)\cong 
  (F^{-1}(0), 0)$} for some holomorphic map germ 
\mbox{$F:(\C^2\!\times T, 0)\to (\C,0)$} with  
$$
F(x,y,\bs)=f(x,y)+\sum\limits^N_{i=1} s_i g_i(x,y,\bs)\,,
$$
where \mbox{$(T,0)$} is a closed subspace of some \mbox{$(\C^N, 0)$}
and \mbox{$\bs=(s_1, \ldots, s_N)$} are 
coordinates of \mbox{$(\C^N, 0)$}. Moreover, under this isomorphism, $\phi$ coincides
with the second projection. We also say that $(i, \phi)$ is isomorphic to the
{\em embedded deformation\/} defined by $F$. A given section
\mbox{$\sigma:(T,0)\to (\sC, 0)$} can always be trivialized, that
is, the ideal \mbox{$I_\sigma= \ker(\sigma^\#\!: 
  \ko_{\sC, 0}\to \ko_{T,0})$} of $\sigma(T,0)$ can be mapped to \mbox{$\langle
x,y\rangle\subset \ko_{\C^2\times T,0}$} under an isomorphism of
embedded deformations. 

The category of deformations (resp.\ of deformations with section) of
$(C,0)$ is denoted by $\Def_{(C,0)}$ (resp.\ by $\Defsec_{(C,0)}$).
The set of isomorphism classes of deformations with section
(over the same base $(T,0)$) is denoted by $\uDefsec_{(C,0)}$
($\uDefsec_{(C,0)}(T,0)$). Here, each isomorphism has to satisfy
\mbox{$\chi=\id_{(T,0)}$}.

\begin{defn}\label{def:DefES}  Let $(C,0)\subset(\C^2, 0)$ be a reduced plane cure
singularity given by $f$ and let $(i, \phi, \sigma)$ be an (embedded) deformation with
section of $(C, 0)$ over $(T,0)$ given by $F$.
The deformation $(i, \phi, \sigma)$ is called 
\begin{enumerate}
\item[$\bullet$] {\em equimultiple (along $\sigma$)\/} if \mbox{$F\in
(I_\sigma)^n$} where \mbox{$n=\ord(f)$} is the multiplicity of $f$
(if $\sigma$ is the trivial section, this means that
\mbox{$\ord_{(x,y)}F=\ord f$}). 
\item[$\bullet$] {\em equisingular (along $\sigma$)\/} if it is equimultiple
along $\sigma$ and if, after blowing up $\sigma$, there exist sections through
the infinitely near points in the first neighbourhood of $(C,0)$ such
that the respective reduced total transforms of $(\sC,0)$ are
equisingular along these sections.  
\end{enumerate}
Further, a deformation of a nodal singularity (with local equation \mbox{$xy=0$}) is
called {\em equisingular\/} if it is equimultiple. (The same applies to
a deformation of a smooth germ.)
\end{defn}

\noindent
Thus, equisingularity of a deformation with section of $(C,0)$ is defined
by induction on the number of blowing ups needed to get a reduced total
transform of $(C,0)$ which consists of nodal singularities only. A deformation without
section is called {\em equisingular\/}, if it is equisingular along some section. 

Let $\Defes_{(C,0)}$, resp.\ $\DefES_{(C,0)}$, denote the category of
equisingular deformations of 
$(C,0)$ as a full subcategory of $\Defsec_{(C,0)}$, resp.\ of $\Def_{(C,0)}$.
 The set of isomorphism classes of equisingular deformations with
 section of $(C,0)$ over 
 $(T,0)$ is denoted by \mbox{$\uDefes_{(C,0)}(T,0)$} and
$$ \uDefes_{(C,0)}\colon
\text{({\it complex germs})}\longrightarrow 
      \Sets\,,\quad (T,0)\longmapsto \uDefes_{(C,0)}(T,0)$$  
is called the {\em functor of equisingular deformations with
  sections}. Similarly, we define $\uDefES_{(C,0)}$, the (abstract)
{\em equisingular deformation functor}.

\noindent
Next, we define deformations of the parametrization.
We fix a commutative diagram of complex (multi-) germs
$$\xymatrix@C=24pt@R=15pt@M=4pt{
(\overline{C},\overline{0})
\ar@{->>}[d]_{n\:} \ar[dr]^-{\varphi}\\
(C,0)\ar@{^{(}->}[r]_-{j} &
(\C^2\!,0)
}
$$
where \mbox{$(C,0)$} is a reduced plane curve
singularity (with a fixed embedding \mbox{$j:(C,0)\hookrightarrow(\C^2,0)$}),
$n$ is its normalization, and \mbox{$\varphi=j\circ
  n$} is its parametrization. If 
\mbox{$(C,0)=(C_1,0)\cup \ldots \cup (C_r,0)$} is the decomposition of
$(C,0)$ into irreducible components, then
\mbox{$(\overline{C},\overline{0})=(\overline{C}_1,\overline{0}_1)\amalg \ldots
  \amalg (\overline{C}_r,\overline{0}_r)$}  is a multigerm, and $n$
maps \mbox{$(\overline{C}_i,\overline{0}_i) \cong (\C,0)$} 
surjectively onto $(C_i,0)$. In particular, by restriction, $n$ induces the
normalization of the component $(C_i,0)$. 

Since \mbox{$(\overline{C},\overline{0})$} and \mbox{$(\C^2\!,0)$} are smooth
(multi-)germs, each deformation of these germs is trivial. 

\begin{defn}\label{def:def of param}
(1) A {\em deformation of the parametrization\/}
\mbox{$(\overline{C},\overline{0})\xrightarrow{\varphi} (\C^2\!,0)$} 
over a germ $(T,0)$ (with compatible sections) is given by the left
(Cartesian) part of the following diagram
\begin{equation}\label{eqn:**}
\xymatrix@M=5pt@C=16pt@R=20pt{ 
(\overline{C},\overline{0}) \ar@{^{(}->}[r]^-{i}\ar@{}[dr]|{\Box} 
\ar[d]_-{\varphi} & 
(\overline{\sC},\overline{0})\ar[d]^-{\phi}  \ar[r]^-{\cong} & 
(\overline{C}\times T,\overline{0}) \ar[d] &
\ar@{=}[l]\coprod\limits_{i=1}^r (\overline{C}_i\times T,\overline{0}_i)\\ 
(\C^2\!,0) \ar@{^{(}->}[r]^-{j}\ar[d]\ar@{}[dr]|{\Box} & 
(\sM,0 )\ar[d]^-{\phi_0}  \ar[r]^-{\cong} & 
(\C^2\!\times T,0) \ar[d]^-{\pr} \\
\{0\} \ar@{^{(}->}[r] & (T,0) \ar@{=}[r] & (T,0)
\ar@/_1pc/@<-1ex>[u]_-{\sigma}
\ar@/_3pc/@<-3ex>[uu]_-{\overline{\sigma}} 
}
\end{equation}

\noindent
where \mbox{$\phi_0\circ\phi$} is flat. We have \mbox{$(\overline{\sC},\overline{0})=
\coprod_{i=1}^r (\overline{\sC}_i,\overline{0}_i)$}, and there are isomorphisms
\mbox{$(\overline{\sC}_i,\overline{0}_i)\cong 
(\overline{C}_i\times T, \overline{0}_i)$}, such that the obvious
diagram (with $\pr$ the projection) commutes.

Systems of compatible  sections
$(\overline{\sigma},\sigma)$ consist of disjoint sections 
{$\overline{\sigma}_i:(T,0)\to (\overline{\sC}_i,\overline{0}_i)$} of
\mbox{$\pr\circ\,\phi_i$} (where \mbox{$\phi_i: 
(\overline{\sC}_i,\overline{0}_i)\to (\sM,0)$} denotes the
restriction of $\phi$) and a section $\sigma$ of $\pr$ such that
\mbox{$\phi\circ\overline{\sigma}_i=\sigma$}, 
\mbox{$i=1,\dots,r$}. 
Morphisms of such deformations  are given by
morphisms of the diagram \eqref{eqn:**}.

\medskip\noindent(2) 
The category of deformations of the
  parametrization $\varphi$ over $(T,0)$ (without sections)
  is denoted by \mbox{$\Def_{(\overline{C},\overline{0}) \to 
    (\C^2\!,0)}(T,0)$}. Its objects are denoted by
$(i,j,\phi,\phi_0)$ or just by $\phi$. 
The corresponding category of deformations of
\mbox{$\varphi$} with
compatible sections is denoted by
\mbox{$\Defsecii_{(\overline{C},\overline{0})\to
      (\C^2\!,0)}(T,0)$}. Its objects are denoted by
  \mbox{$(\phi,\overline{\sigma},\sigma)$}. 
The respective sets of isomorphism classes of deformations are denoted
by \mbox{$\uDef_{(\overline{C},\overline{0}) \to  
    (\C^2\!,0)}(T,0)$} and \mbox{$\uDefsecii_{(\overline{C},\overline{0})\to
      (\C^2\!,0)}(T,0)$}.

\medskip\noindent(3) 
\mbox{$\Tisecii_{(\overline{C},\overline{0})\to 
    (\C^2\!,0)}= \uDefsec_{(\overline{C},\overline{0})\to
    (\C^2,0)}(T_\varepsilon)$} denotes the corresponding vector space
of {\em (first order) 
  infinitesimal deformations\/} of the parametrization with section.
\end{defn}

\noindent
The following theorem shows that deformations of the parametrization induce
(unique) deformations of the equation:

\begin{thm}\label{thm:fitting}
 Each deformation \mbox{$\phi: (\overline{\sC}, 0)\xrightarrow{\pi}
(\C^2\!\times T,0)\xrightarrow{\pr} (T,0)$} of the para\-metri\-zation of the reduced
curve singularity $(C,0)$ induces a deformation of the equation which is
unique up to isomorphism and which is given as follows: the Fitting
ideal of $\pi_\ast\ko_{\Ovkc, \OvO}$, generated by the maximal minors 
of a presentation matrix of $\pi_\ast\ko_{\Ovkc, \OvO}$ as $\ko_{\C^2\times
    T,0}$-module, is a principal ideal which coincides with the kernel
  of the induced morphism of 
  rings \mbox{$\ko_{\C^2\times T,0}\to
  \pi_\ast\ko_{\Ovkc, \OvO}$}. If $F$ is a generator for this ideal, then
$F$ defines an embedded deformation of $(C,0)$. 

In the same way, a deformation
  $(\phi, \Ovsigma, \sigma)$ with compatible sections induces a deformation
  with section of the equation.
\end{thm}

\noindent
The proof of this theorem uses the local criterion of flatness from local
algebra and proceeds by reduction to the special fibre, that is, to
the case that $(T, 0)$ is the reduced point.

A deformation \mbox{$\phi:
  (\overline{C}\times T, 0)\to (\C^2\times T,0)$} of the
parame\-tri\-zation (as in the right-hand part of the diagram
\eqref{eqn:**}) is given by \mbox{$\phi=\{\phi_i=(X_i, Y_i)\}^r_{i=1}$},
$$
\begin{array}{lcl}
X_i (t_i, \bs) & = & x_i(t_i)+A_i(t_i, \bs),\\
Y_i (t_i, \bs) & = & y_i(t_i)+B_i(t_i, \bs),
\end{array}
$$
where \mbox{$X_i, Y_i\in \ko_{\overline{C}\times T, 0}$}, \mbox{$A_i(t_i,
  \bo)\!\!\:=\!\!\:B_i(t_i, \bo)\!\!\: =\!\!\: 0$}, \mbox{$\bs\in
  (T,0)\subset (\C^k\!, 0)$}, 
and where {$\varphi=\{\varphi_i=(x_i, 
  y_i)\}^r_{i=1}$} is the given parametrization of $(C,0)$. We may assume that
    the compatible (multi-)sections
    \mbox{$\overline{\sigma}=\{\Ovsigma_i\}^r_{i=1}$} and $\sigma$ are
    {\em trivial\/}, that is, \mbox{$\Ovsigma_i(\bs)=(\overline{0}_i,
      \bs)$}, \mbox{$\sigma (\bs)=(0,\bs)$}.

\begin{defn} Let $(\phi, \overline{\sigma}, \sigma)\in
\Defsec_{(\overline{C}, \overline{0})\to(\C^2, 0)}(T,0)$ be a deformation of
the parametrization $\varphi: (\overline{C}, 0)\to (\C^2,0)$ as above (with trivial
sections $\Ovsigma, \sigma$).

\begin{enumerate}
\item [(1)] $(\phi, \Ovsigma, \sigma)$ is called {\em equimultiple\/} (along
  $\Ovsigma, \sigma$) if 
$$
\underbrace{\min\{\ord_{t_i}x_i, \ord_{t_i}y_i\}}_{\displaystyle
  =:\ord_{t_i}\varphi_i} 
= \underbrace{\min\{\ord_{t_i}X_i, \ord_{t_i}Y_i\}}_{\displaystyle
  =:\ord_{t_i}\phi_i}\,,\quad i=1, \ldots, r\,. 
$$

\item [(2)] $(\phi, \Ovsigma, \sigma)$ is called {\em equisingular\/} if it
  is equimultiple and if for each infinitely near point
  $p$ of $0$ on the strict transform of $(C,0)$ (after finitely
  many blowing ups) the deformation
  $(\phi, \Ovsigma, \sigma)$ can be lifted to an equimultiple deformation of
  the parametrization of the strict transform in a compatible way (see
  \cite{GLS} for a detailed description of the compatibility condition).
\end{enumerate}
\end{defn}

\noindent
We denote by \mbox{$\Defes_{(\overline{C},\overline{0}) \to
    (\C^2\!,0)}$} the category of equisingular deformations of the parametrization
\mbox{$\varphi:(\overline{C},\overline{0}) \to (\C^2\!,0)$}, and by 
\mbox{$\uDefes_{(\overline{C},\overline{0}) \to (\C^2\!,0)}$} the corresponding
functor of isomorphism classes. Moreover, we
introduce
$$\Ties_{(\overline{C},\overline{0})\to
    (\C^2\!,0)}:=\uDefes_{(\overline{C},\overline{0}) \to
    (\C^2\!,0)}(T_\varepsilon)\,,$$  
the tangent space to this functor.

Note that
  \mbox{$\varphi=(\varphi_i)_{i=1}^r$},
  \mbox{$\varphi_i(t_i)=\bigl(x_i(t_i),y_i(t_i)\bigr)$}, and we set
  $$\dot{\varphi}=
\left\lgroup
\begin{matrix}
\frac{\partial x_1}{\partial t_1}\\[-0.4em]
\vdots\\[-0.7em]
\frac{\partial x_r}{\partial t_r}
\end{matrix}
\right\rgroup
\!\cdot
    \frac{\partial}{\partial x} 
    +\left\lgroup
\begin{matrix}
\frac{\partial y_1}{\partial t_1}\\[-0.4em]
\vdots\\[-0.7em]
\frac{\partial y_r}{\partial t_r}
\end{matrix}
\right\rgroup
\!\cdot\frac{\partial}{\partial y}.$$    
    
\begin{lem}
    With the above notations, there is an isomorphism of vector spaces,
$$ \Ties_{(\overline{C},\overline{0})\to (\C^2\!,0)} \cong \Iesphi \left/
\left(
\dot{\varphi}\cdot\fm_{\overline{C},\overline{0}}+\varphi^\sharp
  (\fm_{\C^2\!,0}) 
\frac{\partial}{\partial x} + \varphi^\sharp (\fm_{\C^2\!,0})
\frac{\partial}{\partial y}
\right)
\right.\,,
$$
where \mbox{$\Iesphi:=\Ies_{(\overline{C},\overline{0})\to
    (\C^2\!,0)}$}
denotes the set of all elements  
$$
\left\lgroup
\begin{matrix}
a_1\\[-0.5em]
\vdots\\[-0.8em]
a_r
\end{matrix}
\!\right\rgroup \!\cdot \frac{\partial}{\partial x}
+
\left\lgroup
\begin{matrix}
b_1\\[-0.5em]
\vdots\\[-0.8em]
b_r
\end{matrix}
\!\right\rgroup \!\cdot \frac{\partial}{\partial y}
\in \fm_{\overline{C},\overline{0}}  \cdot \frac{\partial}{\partial x}
+\fm_{\overline{C},\overline{0}}  \cdot \frac{\partial}{\partial y} $$
such that \mbox{$\bigl\{\bigr(x_i(t_i)+\varepsilon a_i(t_i),
  y_i(t_i)+\varepsilon b_i(t_i)\bigr)\,\big|\, i=1,\dots,r \bigl\}$} 
defines an equisingular deformation of 
\mbox{$\varphi:(\overline{C},\overline{0}) \to (\C^2\!,0)$} over $T_{\varepsilon}$
along the trivial sections.
\end{lem}

\noindent
We call \mbox{$\Iesphi$} the {\em
  equisingularity module of the parametrization\/} of $(C,0)$.
It is an $\ko_{C, 0}$-submodule
of \mbox{$\varphi^\ast \Theta_{\C^2\!,0}= \ko_{\overline{C},\overline{0}}
  \frac{\partial}{\partial x} + 
\ko_{\overline{C},\overline{0}}   \frac{\partial}{\partial y}$}. Here,
$\Theta_{\C^2,0}$ denotes the module of $\C$-derivations
$\Der_\C(\ko_{\C^2,0},\ko_{\C^2,0})$.

The following theorem shows that \mbox{$\uDefes_{(\overline{C},\overline{0})
    \to (\C^2\!,0)}$} is a ``linear'' subfunctor of
\mbox{$\uDefsecii_{(\overline{C},\overline{0})
    \to (\C^2\!,0)}$}. As such, it is already completely determined by its
tangent space. We use the notation
$$ \ba^j = \left\lgroup
\begin{matrix}
a_1^j\\[-0.5em]
\vdots\\[-0.8em]
a_r^j
\end{matrix}
\right\rgroup\,,
\quad 
\bb^j = \left\lgroup
\begin{matrix}
b_1^j\\[-0.5em]
\vdots\\[-0.8em]
b_r^j
\end{matrix}
\right\rgroup\,,
\quad  j=1,\dots,N\,.$$ 

\begin{thm}\label{thm:Defes} With the above notations, the following holds:

\medskip\noindent
(1)\, Let 
\mbox{$(\phi,\osigma,\sigma)$} be a deformation of $\varphi$ with trivial
sections over \mbox{$(\C^N\!,0)$}, where
\mbox{$\phi=\{(X_i,Y_i,\bs)\mid i=1,\dots,r\}$} with 
$$\begin{array}{rcl}
 X_i(t_i,\bs)&=&x_i(t_i)+\displaystyle\sum\limits_{j=1}^N a^j_i(t_i)s_j\,, \quad a^j_i\in
 t_i\C\{t_i\}\,,\\   
 Y_i(t_i,\bs)&=&y_i(t_i)+\displaystyle\sum\limits_{j=1}^N
 b^j_i(t_i)s_j\,, \quad b^j_i\in t_i\C\{t_i\}\,,
\end{array}
$$
\mbox{$i=1,\dots,r$}. Then $\phi$ is equisingular iff
\mbox{$\ba^j \frac{\partial}{\partial x}+\bb^j
  \frac{\partial}{\partial y}\in 
  \Iesphi$} for 
all   \mbox{$j=1,\dots,N$}.

\medskip\noindent
(2)\, Let \mbox{$(\phi,\osigma,\sigma)$} be an equisingular
deformation of $\varphi$ with trivial 
sections over \mbox{$(\C^N\!,0)$}, where
\mbox{$\phi=\{(X_i,Y_i,\bs)\mid i=1,\dots,r\}$} for some \mbox{$X_i,Y_i\in
  \ko_{\C^N\!,0}\{t_i\}$}. Then $(\phi,\osigma,\sigma)$ is a versal
(respectively semiuniversal) object of
\mbox{$\Defes_{(\overline{C},\overline{0}) \to 
    (\C^2\!,0)}$} iff the derivations
$$
\left\lgroup
\begin{matrix}
\frac{\partial X_1}{\partial s_j}(t_1,\bo)\\[-0.4em]
\vdots\\[-0.7em]
\frac{\partial X_r}{\partial s_j}(t_r,\bo)
\end{matrix}
\!\right\rgroup
\cdot \frac{\partial}{\partial x} 
+
\left\lgroup
\begin{matrix}
\frac{\partial Y_1}{\partial s_j}(t_1,\bo)\\[-0.4em]
\vdots\\[-0.7em]
\frac{\partial Y_r}{\partial s_j}(t_r,\bo)
\end{matrix}
\!\right\rgroup
\cdot \frac{\partial}{\partial y} \,, \quad j=1,\dots,N\,,
$$
represent a system of generators (respectively a basis) of the complex
vector space\/
\mbox{$\Ties_{(\overline{C},\overline{0})\to
    (\C^2\!,0)} $}.

\medskip\noindent
(3)\, Let \mbox{$\ba^j \frac{\partial}{\partial x}+\bb^j
  \frac{\partial}{\partial y}\in 
  \Iesphi$}, 
\mbox{$j=1,\dots,N$}, represent a basis (respectively a system of generators)
of\/ \mbox{$\Ties_{(\overline{C},\overline{0})\to
    (\C^2\!,0)}$}. Moreover, let \mbox{$\phi=\{(X_i,Y_i,\bs)\mid
  i=1,\dots,r\}$} be the deformation of $\varphi$ over
\mbox{$(\C^N\!,0)$} given by 
$$
\begin{array}{rcl}
 X_i(t_i,\bs)&=&x_i(t_i)+\displaystyle\sum\limits_{j=1}^N a_i^j(t_i)s_j\,, \\   
 Y_i(t_i,\bs)&=&y_i(t_i)+\displaystyle\sum\limits_{j=1}^N b_i^j(t_i)s_j\,, 
\end{array}
$$
\mbox{$i=1,\dots,r$}, and let $\osigma,\sigma$ be the trivial
sections. Then $(\phi,\osigma,\sigma)$ is a 
semiuniversal (respectively versal) 
equisingular deformation of $\varphi$ over \mbox{$(\C^N\!,0)$}. 
In particular, equisingular deformations of the
parametrization are unobstructed, and the semiuniversal deformation has a
smooth base space of dimension 
\mbox{$\dim_\C \Ties_{(\overline{C},\overline{0})\to
    (\C^2\!,0)}$}.
\end{thm}

\noindent
In the proof we make a power series ''Ansatz'' and then we explicitly verify the condition of versality in the spirit of Schlessinger.

To compute a semiuniversal equisingular deformation of
\mbox{$\varphi:(\overline{C},\overline{0}) \to (\C^2\!,0)$}, we only
need to compute a basis of \mbox{$\Ties_\varphi$} by Theorem
\ref{thm:Defes}. Moreover, if all branches of $(C,0)$ have different
tangents, then 
$\Ties_\varphi$ decomposes as 
$$   \Ties_\varphi = \bigoplus_{i=1}^r \Ties_{\varphi_i}\,,$$
where $\varphi_i$ is the parametrization of the $i$-th branch of $(C,0)$. In
general, \mbox{$\Ties_\varphi$} can be computed following the lines of the
proof of Theorem \ref{thm:Defes}.

\begin{examples}
  (1) Consider the parametrization \mbox{$\varphi:t\mapsto (t^2,t^7)$}
    of an $A_6$-sin\-gu\-la\-rity. A basis for the mo\-dule of
    equimultiple deformations $M^{em}_\varphi$ is given by
    \mbox{$\{ t^3\frac{\partial}{\partial 
  y},t^5\frac{\partial}{\partial y} \}$}. 
Blowing up the trivial section of the deformation of $\phi$ given by
\mbox{$X(t,\bs)=t^2$}, \mbox{$Y(t,\bs)=t^7\!+s_1 t^3\!+s_2 t^5$}, we 
get 
$$U(t,\bs)=t^2\,,\quad V(t,\bs)=\frac{Y(t,\bs)}{X(t,\bs)}=t^5\!+s_1 t+s_2
  t^3\,,$$ 
which is equimultiple along the trivial section iff
\mbox{$s_1=0$}. Blowing up once more, we get the necessary condition
\mbox{$s_2=0$} for equisingularity. Hence, \mbox{$\Ties_{\varphi}=0$}
as expected for a simple singularity (each equisingular deformation of
a simple singularity is known to be trivial).  

\medskip\noindent
(2)\,  For the parametrization \mbox{$\varphi:t\mapsto (t^3,t^7)$}
of an $E_{12}$-singularity, a basis for \mbox{$M^{\text{\it em}}_\varphi$} is given by
\mbox{$\{ t^4\frac{\partial}{\partial 
  y},t^5\frac{\partial}{\partial y},t^8\frac{\partial}{\partial y}
\}$} (resp.\ by \mbox{$\{ t^4\frac{\partial}{\partial 
  x},t^4\frac{\partial}{\partial  y},t^5\frac{\partial}{\partial y}
\}$}). Blowing up the trivial section, only 
$t^8\frac{\partial}{\partial y}$ (resp.\  $t^4\frac{\partial}{\partial 
x}$) survives for an equimultiple deformation. 
It also survives in further blowing ups. Hence, \mbox{$X(t,s)=t^3$},
\mbox{$Y(t,s)=t^5+st^8$} (resp.\  \mbox{$X(t,s)=t^3+st^4$},
\mbox{$Y(t,s)=t^5$}) is a semiuniversal equisingular deformation of
$\varphi$.
\end{examples}

\noindent
The following theorem relates deformations of the parametrization to the
$\delta$-constant stratum in the semiuniversal deformation of the
equation. It is an improvement of the results by Teissier and Raynaud,
by Chian-Hsieh and Lipman, and by Diaz and Harris
\cite{DiH}. 

\begin{thm} With the above notations, the following holds:

\medskip\noindent
(1)\, Let \mbox{$(\overline{\sC}, \overline{0})\to (\sM,0)\to (T,0)$} be
  a deformation of \mbox{$\varphi:(\overline{C}, \overline{0})\to
    (\C^2,0)$}, and let \mbox{$(\sC, 
  0)=\phi(\overline{\sC}, \overline{0})\to (T,0)$} be the induced 
deformation of the equation of $(C,0)$. Then \mbox{$\delta(\sC_t)=\sum_{x\in
  \Sing (\sC_t)} \delta(\sC_t,x)$} is constant for \mbox{$t\in T$} near
$0$.\footnote{For germs $(\sC,0), (T,0)$, etc., $\sC,T$, etc. always
  denote sufficiently small representatives.} 

\medskip\noindent
(2)\, Let \mbox{$(\sC,0)\to (T,0)$} be the semiuniversal deformation of the
    equation of $(C,0)$, and let \mbox{$\Delta^\delta: = \{t\in T\ |\
      \delta(\sC_t)=\delta(C,0)\}$}
be the {\em $\delta$-constant stratum\/} of $(C,0)$. Then:
\begin{enumerate}
\item [(a)]  The semiuniversal deformation of the parame\-tri\-zation of
  $(C,0)$ is induced from \mbox{$(\sC,0)\to (T,0)$} via a morphism
  {$\Phi:(S,0)\to (T,0)$} such  
  that 
\begin{quote}
\begin{itemize}
\item \mbox{$\Phi(S,0)=(\Delta^\delta, 0)$} and
\item \mbox{$\Phi:(S,0)\to (\Delta^\delta, 0)$} is the normalization of
  $(\Delta^\delta, 0)$.
\end{itemize}
\end{quote}
\item [(b)] $(\Delta^\delta, 0)$ has a smooth normalization and $\codim_{(T,0)}
  (\Delta^\delta, 0)=\delta$.
\item [(c)] $(\Delta^\delta, 0)$ is smooth iff all branches $(C_i, 0)$ of $(C,0)$
  are smooth.
\end{enumerate}
\end{thm}

\noindent
The proof of this theorem uses the results of \cite{Tei1} and \cite{CHL} 
  mentioned in the introduction, the fact that for every plane curve
  singularity $(C,0)$ there is a $\delta$-constant deformation such that 
  the general fibre has $\delta(C,0)$ simple nodes, and an exact sequence
  relating first order deformations of the equation and of the parametrization.

When passing to equisingular deformations of the parametrization, we
have to consider deformations with compatible sections. It can be
shown that the sections are unique (in characteristic $0$). Then, a
refinement of the above arguments for equisingular deformations proves the 
following theorem:

\begin{thm}\label{thm:esdef param, eqn} 
Let \mbox{$(\overline{\sC}, \overline{0})\to (\sM,0)\to (S,0)$}
be the semiuniversal equisingular deformation of the parametrization of
$(C,0)$, and let \mbox{$\Phi:(S,0)\to (T,0)$} be the inducing morphism
to the base space of the semiuniversal deformation of the equation.
Then $\Phi$ is an isomorphism onto the $\mu$-constant
stratum \mbox{$(\Delta^\mu,0)\subset (T,0)$}. In particular, $(\Delta^\mu,0)$ is
smooth.
\end{thm}

\section{The Algorithms}
\setcounter{thm}{0}
\setcounter{equation}{0}

\noindent
The idea of the following algorithm to compute the equisingularity stratum of
a family of plane curve singularities with trivial section was developed in our joint
preprint \cite{CGL}. In that paper we introduced the notion of equisingularity
for plane algebroid curves given by a formal power series \mbox{$f\in K[[x,y]]$},
where $K$ is an algebraically closed field of any characteristic.

The definitions of the previous section remain true, mutatis
mutandis, for algebroid curves. However, we cannot use the geometric
language. Instead of morphisms 
between complex space germs, we have to consider morphisms (in the
opposite direction) between the corresponding local analytic algebras.
Points \mbox{$t\in T$} close to $0$ have to be 
replaced by generic points of \mbox{$\Spec\ko_{T,0}$}, etc. For
\mbox{$K=\C$}, it does not make any difference whether we consider
convergent or formal power series. The 
reason for considering convergent power series in the previous section is
that the concept of equisingularity can be best explained in a geometric
context and that a great deal of the motivaton comes from topology.

However, there is an important difference between the case of
characteristic $0$ and the case of positive characteristic. As 
shown in \cite{CGL},  in positive 
characterestic we have two equally important notions of
equisingularity, namely {\em weak\/}
and {\em strong equisingularity\/} which coincide in characteristic $0$. The
definitions for equisingularity given in Section \ref{sec:theorems}
(appropriately formulated on the level of analytic rings), either for
the equation or for the parametrization, refer to the notion of strong
equisingularity (which we continue to call equisingularity).

The theorems of the previous section remain true for algebraically closed
fields $K$ of characteristic $p$ as long as $p$ does not divide the multiplicity of
any factor of \mbox{$f\in K[[x,y]]$} (in particular, for each
algebraically closed field of characteristic $0$). This
result, proved in \cite{CGL1} has the important computational consequence that
for a power series $f$ with integer coefficients we can compute characteristic
numerical invariants like $\delta, r$, and the Puiseux
pairs\footnote{Note that, in  positive characteristic, the Milnor
  number as defined on Page \pageref{page:Milnor} de\-pends on the equation $f$ 
  and not only on the ideal $\langle f\rangle$. Instead, we define the
  Mil\-nor number in characteristic $p$ as
  \mbox{$\mu:=2\delta-r+1$}.}in characteristic $0$ by 
computing them modulo a prime number $p$, where $p$ is bigger than the
multiplicity of $f$. This is the reason why we work in this section with
analytic local rings over a field $K$ of possibly positive
characteristic. 

In \cite{CGL1}, we treat the case of arbitrary characteristic. Here,
we treat only (strong) equisingularity and assume, 
that the characteristic of $K$ does not divide the multiplicity of any
branch of $(C,0)$. 

Since the Puiseux expansion is in general not available in positive
characteristic, we work with the Hamburger-Noether expansion instead
(cf. \cite{Cam,Cam1}).

\smallskip\noindent
We fix the notations. $K$ denotes an algebraically closed field of characteristic \mbox{$p \ge 0$}.
All rings in this section will be Noetherian complete local $K$-algebras $A$
with maximal ideal $\fm_A$ such that \mbox{$A/\fm_A = K$}.  The category of
these algebras is denoted by $\sA_K$. Further, we denote by $K[\varepsilon]$
the two-dimensional $K$-algebra with \mbox{$\varepsilon^2 = 0$}.
Let $C$ be a reduced algebroid plane curve singularity over $K$, defined
by the (square-free) power series \mbox{$f \in K[[x,y]]$}.
\[
R = \ko_C = P/\langle f \rangle,\quad P = K[[x,y]],
\]
denotes the complete local ring of $C$. Let \mbox{$f=f_1\cdot\ldots\cdot f_r$}
be an irreducible factorization of $f$. The rings 
$$
R_i = P/\langle f_i\rangle\,,\quad i=1,\dots,r\,,
$$
are the complete local rings of the branches $C_i$ of $C$.
The normalization $\Rbar$ of $R$ is the integral closure of $R$ in its
total ring of fractions $\Quot(R)$. It is the direct sum of the
normalizations $\Rbar_i$ of $R_i$, \mbox{$i=1,\dots,r$}, hence a semilocal
ring. Each $\Rbar_i$ is a discrete valuation ring, and we can choose
uniformizing parameters $t_i$ such that \mbox{$\Rbar_i\cong K [[t_i]]$}. After
fixing the parameters $t_i$, we identify $\Rbar_i$ with \mbox{$K [[t_i]]$} and 
get 
$$
\Rbar = \bigoplus^r_{i=1} \Rbar_i =  \bigoplus^r_{i=1} K [[t_i]]\,.
$$
The normalization map \mbox{$R \hookrightarrow \overline{R}$} is induced by
a mapping \mbox{$\varphi:P \to \overline{R}$},\, {$(x,y) \mapsto (x_i(t_i),
  y_i(t_i))^r_{i=1}$}, which is called a {\em parametrization} of $R$.

The following definition is to local analytic $K$-algebras what Definition
\ref{def:def of param} is to analytic germs:

\begin{defn}\label{def:1.1}
  A {\em deformation with sections of the parametrization\/} of $R$ over
  \mbox{$A\in \sA_K$} is a commutative diagram with Cartesian squares  
$$
\UseComputerModernTips
\xymatrix@C=9pt@R=2pt@M=6pt{
{\Rbar}\ar@{}[ddrr]|{\Box} && \ar@{->>}[ll]
{\Rbar}_A\ar@/^3pc/[dddd]^{\osigma=\{\osigma_i\mid i=1,\dots,r
  \}}\\   
\\
P \ar[uu]^{\varphi}\ar@{}[ddrr]|{\Box} && \ar@{->>}[ll] P_A
\ar[uu]_{\varphi_A} \ar@/^/[dd]^{\sigma} \\ 
 \\
 K \ar@{^{(}->}[uu] && \ar@{->>}[ll] A
\ar[uu]
}
$$
with \mbox{$\Rbar_A=\bigoplus_{i=1}^r \Rbar_{A,i}$}, where $\Rbar_{A,i}$,
\mbox{$i=1,\dots,r$}, and $P_A$ are Noetherian complete local $K$-algebras
which are flat over $A$. $\sigma$ is a section of \mbox{$A\to P_A$}, and
$\osigma_i$ is a section of \mbox{$A\to \Rbar_{A,i}$}, \mbox{$i=1,\dots,r$}. We
denote such a deformation by \mbox{$\xi=(\varphi_A,\osigma,\sigma)$}.  

A morphism from \mbox{$\xi$} to another deformation
\mbox{$(P_B\!\xrightarrow{\varphi_B}\! {\Rbar}\!\,'_B,\osigma_B,\sigma_B)$}
over \mbox{$B\in \sA_K$} is then given by morphisms 
\mbox{$A\to B$}, \mbox{$P_A\to P_B$} 
and \mbox{$\Rbar_{A,i}\to \Rbar_{B,i}$} in $\sA_K$ such that the resulting
diagram commutes. The category of such deformations is denoted by
$\Defsec_{\Rbar\leftarrow P}$. If we consider only deformations over a fixed
base $A$, we obtain the (non-full) subcategory \mbox{$\Defsec_{\Rbar\leftarrow
    P}(A)$} with morphisms being the identity on $A$.
\mbox{$\Defsec_{\Rbar\leftarrow P}$} is a fibred gruppoid over $\sA_K$, that
is, each morphism in  \mbox{$\Defsec_{\Rbar\leftarrow
    P}(A)$} is an isomorphism.
\end{defn}

\noindent
Since $P$ and the $\Rbar_i$ are regular local rings, each deformation of $P$
  and of $\Rbar$ is trivial. That
  is, there are isomorphisms \mbox{$P_A\cong A[[x,y]]$}  and
  \mbox{$\Rbar_A\cong \bigoplus_{i=1}^r A[[t_i]]$} over $A$, mapping 
  the sections $\sigma$ and $\osigma_i$ to the trivial sections. Hence, each
  object in $\Defsec_{\Rbar\leftarrow P}(A)$ is isomorphic to a diagram of the
  form  
\begin{equation}\label{eqn:def of param}
\UseComputerModernTips
\xymatrix@C=9pt@R=14pt@M=6pt{
\bigoplus\limits_{i=1}^r K[[t_i]]\ar@{}[drr]|{\Box} && \ar@{->>}[ll]
\bigoplus\limits_{i=1}^r
A[[t_i]]\ar@/^3pc/[dd]^{\osigma=\{\osigma_i\,\mid i=1,\dots,r 
  \}}\\
K[[x,y]] \ar[u]^-{\varphi}\ar@{}[drr]|{\Box} && \ar@{->>}[ll] A[[x,y]]
\ar[u]_-{\varphi_A} \ar@/^/[d]^{\sigma} \\ 
 K \ar@{^{(}->}[u] && \ar@{->>}[ll] A
\ar[u]
}
\end{equation}
where $\varphi_A$ is the identity on $A$ and $\sigma$, $\osigma_i$ are the
trivial sections (that is, the cano\-nical epimorphisms mod $x,y$, respectively
mod $t_i$). Here, 
$\varphi_A$ is given by {$\varphi_A=(\varphi_{A,1},\dots,\varphi_{A,r})$},
where $\varphi_{A,i}$ is determined by 
$$ \varphi_{A,i}(x)=X_i(t_i)\,,\quad \varphi_{A,i}(y)=Y_i(t_i)\in
t_iA[[t_i]]\,,$$ 
\mbox{$i=1,\dots,r$}, such that \mbox{$X_i(t_i)\equiv x_i(t_i)$},
\mbox{$Y_i(t_i)\equiv y_i(t_i)$} mod $\fm_A$. 

We write $\uDefsec_{\Rbar\leftarrow P}(A)$ for the set of isomorphism
classes of objects in $\Defsec_{\Rbar\leftarrow P}(A)$, and we
denote by \mbox{$\uDefsec_{\Rbar\leftarrow P} : \sA_K \to \text{(Sets)}$}
the corresponding deformation functor. Moreover, we denote by
\mbox{$\Tisec_{\Rbar\leftarrow P}:=\uDefsec_{\Rbar\leftarrow P}(\Keps)$} the
vector space of {\em (first order) infinitesimal deformations\/} of the
parametrization of $R$.

\begin{rem}
Replacing in the above definition the parametrization
\mbox{$P\xrightarrow{\varphi}  \Rbar$} by the normalization
\mbox{$R\hookrightarrow \Rbar$}, we get the 
functor $\uDefsec_{\Rbar\leftarrow R}$ of {\em deformations of the
  normalization\/}. The version of Theorem \ref{thm:fitting} for local
$K$-algebras implies that this functor is naturally equivalent to
$\uDefsec_{\Rbar\leftarrow P}$. 
\hfill \qedsymbol
\end{rem}

\noindent
It is now straightforward to translate the definition of equisingular
deformations of the parametrization from the geometric to the algebraic
context. We leave this to the reader. For the algorithms, it is only important
to know that a deformation \eqref{eqn:def of param} is equisingular iff (up
to a reparametrization) it is given by a Hamburger-Noether deformation
of $C$ over $A$, which we introduce next (see Proposition
\ref{prop:ES=HN} below). 

\begin{defn} A {\em Ham\-bur\-ger-Noether expansion (HNE)\/} $\kh_A$ over $A$ 
is a finite system of equations in the variables 
$z_{-1},z_0,\dots, z_{\ell}$ of type

\begin{center}
\hfill \qquad
$
\renewcommand{\arraystretch}{1.2} 
\begin{array}{rcl}
z_{-1}& = & a_{0,1}z_0 + a_{0,2}z_0^2 + \ldots + a_{0,d_0}z_0^{d_0}\!
+z_0^{d_0}z_1\\  
z_{0}& = & \phantom{a_{0,1}z_0 + }\;\, a_{1,2}z_1^2 + \ldots +
a_{1,d_1}z_1^{d_1}\! +z_1^{d_1}z_2\\
\vdots\phantom{i} & & \phantom{a_{0,1}z_0 + A} \vdots \\
 z_{j-1}& = & \phantom{a_{0,1}z_0 + }\;\, a_{j,2}z_j^2 + \ldots +
 a_{j,d_j}z_j^{d_j}\! +z_j^{d_j}z_{j+1} \\
\vdots\phantom{i} & & \phantom{a_{0,1}z_0 + A}\vdots \\
 z_{\ell-2}& = & \phantom{a_{0,1}z_0 + }\;\, a_{\ell-1,2}z_{\ell-1}^2
 + \ldots +  a_{\ell-1,d_{\ell-1}}z_{\ell-1}^{d_{\ell-1}}\!
 +z_{\ell-1}^{d_{\ell-1}}z_{\ell} \\ 
z_{\ell-1} & = &  \phantom{a_{0,1}z_0 + }\;\, a_{\ell,2}z_{\ell}^2 +
a_{\ell,3}z_{\ell}^3 + \ldots \ldots \ldots \ldots \,,
\end{array}
$
\hfill $(\kh_A)$
\end{center}

\noindent
where $\ell$ is a nonnegative integer, the coefficients $a_{j,k}$ are
elements of $A$, the \mbox{$d_j$} are positive integers, and 
we assume that the first nonzero coefficient in each row, except in the
first one, is a unit in $A$. 
Finally, if \mbox{$\ell>0$}, then the power series
\mbox{$H_{A,\ell}(z_{\ell}):=\sum_{k=2}^\infty 
  a_{\ell,k}z_{\ell}^i$} on the right-hand side of the last equation in 
$\kh_A$ is nonzero. We call $\ell$ the {\em length\/} of $\kh_A$.

Given a Hamburger Noether expansion $\kh_A$ over $A$, we
define the {\em residual HNE\/} $\Res(\kh_A)$ to be the Hamburger-Noether
expansion over $K$ 
obtained by substituting the coefficients \mbox{$a_{j,k}\in A$} by the
respective resid\-ual classes \mbox{$(a_{j,k}\:\mod \fm_A)$}. 
\end{defn}

\begin{rem}
Let $C$ be as above, and let  
\[
\Lambda = \Lambda_1 \cup \Lambda_2
\]
be the partition of the index set \mbox{$\Lambda = \{1, \dots, r\}$} such that 
$\Lambda_1$ (resp.\ $\Lambda_2$) consists of those indices $k$ for
which the line \mbox{$\{x = 0\}$} is transversal (resp.\ tangent) to the
branch $C_{i}$. Then associated with each branch $C_{i}$ one has a unique
Hamburger-Noether expansion $\kh^{(i)}_K$ over $K$ of 
some length $\ell_i$ such
that, setting \mbox{$y := z_{-1}$}, \mbox{$x := z_0$} if \mbox{$i \in
  \Lambda_1$}  and \mbox{$x := z_{-1}$},  \mbox{$y :=
z_0$} if \mbox{$i \in \Lambda_2$}, and \mbox{$t:=z_{\ell^{(i)}}$},
and making successive back-substitutions in $\kh^{(i)}_K$,
we obtain power series \mbox{$x(t),y(t)\in K[[t]]$} defining a
parametrization of the branch $C_i$.
The uniqueness comes from the fact that, since a
transversal parameter is fixed, the data of the Hamburger-Noether expansion
$\kh^{(i)}_K$ collect the information about the coordinates of the successive
infinitely near points on the branch $C_i$ in appropriate coordinate systems
(see \cite[Ch.\ II]{Cam}). Further, 
the expansions $\kh^{(i)}_K$ are pairwise different in $\Lambda_1$
and in $\Lambda_2$, and for \mbox{$i \in \Lambda_2$} one has, in addition to
the defining properties for a Hamburger-Noether expansion, that
\mbox{$a^{(i)}_{01}=0$}. 
\phantom{a}\hfill\qedsymbol
\end{rem}

\begin{defn}
 A {\em deformation of the Hamburger-Noether expansion of\/ $C$
  over $A$\/} (or simply a {\em
   Hamburger-Noether deformation of\/ $C$ over $A$}\/) is a system of
 Ham\-bur\-ger-Noether expansions \mbox{$\kh_A^{(i)}$},
 \mbox{$i=1,\dots,r$}, over $A$,

\begin{center}
\hfill \qquad
$
\renewcommand{\arraystretch}{1.6} 
\begin{array}{rcll}
 z_{j-1}& = & H_{A,j}^{(i)}(z_j) + z_j^{d_j^{(i)}}z_{j+1} \,,&
 \quad  j=0,\dots, \ell^{(i)}\!-1\,,\\
 z_{\ell^{(i)}\!\!\;-1}& = & H_{A,\ell^{(i)}}^{(i)}(z_{\ell^{(i)}})\,, 
\end{array}$
\hfill $(\kh_A^{(i)})$
\end{center}

\noindent
such that, for each  \mbox{$i\neq i'\in \{1,\dots,r\}$} the following holds:
\begin{description}
\item[(HN1)] \mbox{$\Res(\kh_A^{(i)})=\kh^{(i)}_K$}, the Hamburger-Noether
  expansion for $C_i$ (over $K$).
\item[(HN2)] If $i$ and $i'$ are either both in $\Lambda_1$ or both in
  $\Lambda_2$ and if $j_0$ denotes the smallest integer such that
  \mbox{$(d^{(i)}_{j_0}, H_{A,{j_0}}^{(i)})\neq (d^{(i')}_{j_0},
    H_{A,{j_0}}^{(i')})$}, then either the multiplicity of
  \mbox{$H_{A,{j_0}}^{(i)}\!- H_{A,{j_0}}^{(i')}\in A[[z_{j_0}]]$}
  exceeds the minimum of $d^{(i)}_{j_0},\,d^{(i')}_{j_0}$, or the
  coefficient of its term of smallest degree is a unit in $A$.
\end{description}
\end{defn}

\begin{exmp}
  Let \mbox{$K=\C$} and \mbox{$A=\C[[s]]$}. Then the system
$$
(\kh_A^{(1)}) \ 
\renewcommand{\arraystretch}{1.2} 
\renewcommand{\arraycolsep}{3pt} 
\begin{array}{rcl}
   z_{-1} &=& sz_0+z_0^2z_1\\
   z_0 &=& z_1z_2 \\
   z_1 &=& (1+s)z_2^3
 \end{array}
\quad
(\kh_A^{(2)}) \ 
\renewcommand{\arraystretch}{1.2} 
\begin{array}{rcl}
   z_{-1} &=& sz_0+z_0^2z_1\\
   z_0 &=& z_1z_2 \\
   z_1 &=& (1+s)z_2^3+z_2^7+s^2z_2^8+\displaystyle\sum\limits_{k=0}^{\infty} z_2^{11+4k}  
 \end{array}
$$
is a Hamburger-Noether deformation of
\mbox{$C=\{(y^4\!-x^{11})(y^4\!-x^{11}\!-x^{12})=0\}$} 
over $A$. If we replace the last equation in $\kh_A^{(1)}$ by
\mbox{$z_1 = z_2^3$}, then $\kh_A^{(1)}$ is still a Hamburger-Noether expansion
over $A$, but $\kh_A^{(1)},\kh_A^{(2)}$ do not define a Ham\-bur\-ger-Noe\-ther
deformation of $C$ over $A$ (the condition (HN2) is not satisfied).
\end{exmp}

\smallskip
\begin{rem}
By setting 
$$ Y_i:=z_{-1}\,, \
X_i:=z_0  \,\text{ for } i\in \Lambda_1, \quad 
X_i:=z_{-1}\,, Y_i:=z_0\, \text{ for } i\in \Lambda_2\,,
$$
and \mbox{$t_i:=z_{\ell^{(i)}}$}, and by making successive back-substitutions,
we obtain power series {$X_i(t_i),Y_i(t_i)\in A[[t_i]]$}, \mbox{$i=1,\dots,r$},
satisfying \mbox{$X_i(0)=Y_i(0)=0$}. These define a deformation of the
parametrization 
$$
\varphi : P\to \Rbar=\bigoplus_{i=1}^r K[[t_i]]\,, \quad 
(x,y) \mapsto \bigl(x_i(t_i),y_i(t_i) \bigr)_{i=1}^r,$$
 \mbox{$x_i(t_i):=X_i(t_i) 
\:\mod \fm_A$}, \mbox{$ y_i(t_i):=Y_i(t_i) \:\mod \fm_A$}, of $C$ which is
induced by the system of Hamburger-Noether expansions
$\kh^{(1)}_K,\dots,\kh^{(r)}_K$ for $C$. 

For instance, in the above example, we get the deformation of the parametrization
given by 
$$ \begin{array}{rcl}
\bigl(X_1(t_1), Y_1(t_1)\bigr) & = & \bigl((1+s)t_1^4,\,
(s+s^2)t_1^4+(1+s)^3t_1^{11} \bigr)\,, \\[0.3em]
\bigl(X_2(t_2), Y_2(t_2)\bigr) & = &
\bigl((1+s)t_2^4+t_2^8+s^2t_2^9+t_2^{12}+\ldots,\,\\[0.1em] 
&& \quad(s+s^2)t_2^4+st_2^8+s^3t_2^9+(1+s)^3t_2^{11}+st_2^{12}+\ldots
\bigr)\,. 
\end{array}
$$

\vspace{-14pt}
 \hfill\qedsymbol
\end{rem}

\smallskip
\begin{prop}\label{prop:ES=HN}
The deformation of the parametrization \mbox{$\varphi:P\to \Rbar$}
associated to a Hamburger-Noether deformation of $C$ over $A$ is equisingular
(along the trivial section $\sigma$). This association is functorial 
in $A$. Conversely, every equisingular deformation of the parametrization with
trivial section $\sigma$ is given, up to a re-parametrization, by a
Hamburger-Noether deformation. 
\end{prop}

\noindent
The proof of this proposition (as given in \cite{CGL}) provides an
algorithm for finding the Hamburger-Noether deformation of $C$
associated to an equisingular deformation of the parametrization. 
This leads to the following algorithm which allows one to decide whether 
a given deformation of the parametrization is equisingular:

\begin{alg}[Check equisingularity]\label{alg:1}\hfill

\begin{entry}
\item[Input] \mbox{$X_i(t_i),Y_i(t_i) \in A[[t_i]]$}, \mbox{$i=1,\dots,r$},
  defining a deformation of the para\-met\-rization of a reduced plane curve
  singularity over a complete local $K$-algebra
  \mbox{$A=K[[s_1,\dots,s_N]]/I$}. 
\item[Output] $1$ if the deformation is equisingular along the trivial
  section, $0$ otherwise.
\end{entry}

\medskip\noindent
{\it Step 1.}\, {\it (Initialization)}

\smallskip
\begin{itemize}
\item For each \mbox{$i=1,\dots,r$}, set
$$x_i(t_i):= (X_i(t_i) \:\mod
    \fm_A)\,, \quad y_i(t_i):= (Y_i(t_i) \:\mod \fm_A)\,.$$
\item Set \mbox{$\Lambda_1:= \{i \mid \ord x_i(t_i)\leq \ord y_i(t_i)\}$},\;
  \mbox{$\Lambda_2:=\{1,\dots,r\}\setminus \Lambda_1$}.
\end{itemize}

\medskip\noindent
{\it Step 2.}\, If for some \mbox{$1\leq i\leq r$} the condition
$$\begin{array}{rcl}
\ord x_i(t_i) =
\ord_{t_i} X_i(t_i) \le \ord_{t_i} Y_i(t_i) &  \text{ if }\: i \in \Lambda_1\,,\\
\ord y_i(t_i) = \ord_{t_i} Y_i(t_i) \le \ord_{t_i} X_i(t_i) & \text{ if
}\: i \in \Lambda_2\,.   
\end{array}
$$
is not fulfilled then {\sc Return}(0).

\medskip\noindent
{\it Step 3.}\, {\it (Compute the Hamburger-Noether expansions 
$\kh^{(1)}_A,\dots,\kh^{(r)}_A$)}

\smallskip
\noindent
For each \mbox{$i=1,\dots,r$} do the following:

\smallskip
\begin{itemize}
\itemsep4pt
\item Set \mbox{$Z_0 := X_i(t_i)$}, \mbox{$Z_{-1}:= Y_i(t_i)$} if
\mbox{$i \in \Lambda_1$}, and 
\mbox{$Z_0 := Y_i(t_i)$}, \mbox{$Z_{-1} := X_i(t_i)$} if \mbox{$i \in
  \Lambda_2$}. 
\item If \mbox{$\ord_{t_i} Z_0=1$}, then the Hamburger-Noether expansion
  $\kh^{(i)}_A$ has length \mbox{$\ell^{(i)}=0$} and the coefficients
$a^{(i)}_{0,k}$ are obtained by expanding $Z_{-1}$ as a 
  power series in $Z_0$. 
\item Set \mbox{$j:=0$}, \mbox{$k:=0$}. 
\item While \mbox{$\ord_{t_i} Z_j>1$}  do the following:

\smallskip
\begin{itemize}
\itemsep3pt
\item While \mbox{$\ord_{t_i} Z_{j-1}\geq \ord_{t_i} Z_{j}$}, set
  \mbox{$k:=k+1$}, define \mbox{$a_{j,k}^{(i)}\in A$} to be the residue
  modulo $t_i$ of $Z_{j-1}/Z_j$, and set
 $$Z_{j-1} := \frac{Z_{j-1}}{Z_j} - a_{j,k}^{(i)}\in A[[t_i]]\,.$$
\item If the leading coefficient of $Z_{j-1}$ is not a unit in $A$, then 
{\sc Return}(0).
\item Set \mbox{$d_j^{(i)}:=k$}, \mbox{$Z_{j+1}:=Z_{j-1}$}, and
  \mbox{$j:=j+1$}.
\end{itemize}
\item The Hamburger-Noether expansion
  $\kh^{(i)}_A$ has length \mbox{$\ell^{(i)}=j$} and the coefficients
$a^{(i)}_{j,k}$ in its last row are obtained by expanding $Z_{j-1}$ as a 
  power series in $Z_j$.
\end{itemize}

\medskip\noindent
{\it Step 4.}\, {\it (Check condition (HN2) for a Hamburger-Noether
    expansion)}

\smallskip\noindent
For each \mbox{$i=1,\dots,r$}, \mbox{$j=1,\dots,\ell^{(i)}$}, set
\mbox{$H_{A,j}^{(i)}:=\sum_{k} 
  a^{(i)}_{j,k}z_j^k\in A[[z_j]]$}. If the condition (HN2) is 
satisfied then {\sc Return}(1), otherwise {\sc Return}(0). \hfill\qedsymbol
\end{alg}

\smallskip
\begin{rem} Algorithm \ref{alg:1} can be extended in an obvious way to an
  algorithm which computes for an arbitrary deformation with trivial
  section of the parametrization 
  of $C$ over $A$ an ideal \mbox{$\fa \subset A$} such that the induced
  deformation over \mbox{$A/\fa$} is 
equisingular and, if \mbox{$\fb \subset A$} is any other ideal with this
property, then \mbox{$\fb \supset \fa$}.
If we apply this algorithm to the deformation of the parametrization
  given by $$X_i(t) = x_i(t_i) + \sum_{k=1}^N \varepsilon_k
  a_i^k(t_i)\,,\quad Y_i(t_i) = y_i(t_i) + \sum_{k=1}^N \varepsilon_k
  b_i^k(t_i)$$ 
over the Artinian $K$-algebra
\mbox{$K[\beps]/\langle \beps\rangle^2$},
\mbox{$\beps=(\varepsilon_1,\dots,\varepsilon_N)$}, where 
  the 
$$(\ba^j,\bb^j) \in  
\bigoplus_{i=1}^r \bigl(t_iK[[t_i]] \oplus t_iK[[t_i]]\bigr)\,,\quad
k=1,\dots,N\,,$$  
represent a $K$-basis of \mbox{$\Tisec_{\Rbar \leftarrow P}$}, then
the conditions obtained are $K$-linear equations in the $\varepsilon_k$.
Solving the system of these linear equations and
  restricting the family to the corresponding subspaces, we get
  a family 
$$\widetilde{X}_i(t) = x_i(t_i) + \sum_{k\in I} \varepsilon_k
  \widetilde{a}\!\,_i^k(t_i)\,,\quad \widetilde{Y}_i(t_i) = y_i(t_i) +
  \sum_{k\in I} \varepsilon_k 
  \widetilde{b}\!\,_i^k(t_i)\,,$$ 
where $I$ is a subset of \mbox{$\{1,\dots,N\}$}, and where the
\mbox{$(\widetilde{\ba}\!\,^k, \widetilde{\bb}\!\,^k)$} are $K$-linear
combinations of the $(\ba^k,\bb^k)$. Then the \mbox{$(\widetilde{\ba}\!\,^k,
  \widetilde{\bb}\!\,^k)$}, \mbox{$k\in I$}, generate a linear subspace $T$ of
\mbox{$\Tisec_{\Rbar \leftarrow P}$} which is necessarily equal to
  \mbox{$\Ties_{\Rbar \leftarrow P}$}.  This follows, since \mbox{$T
    \subset \Ties_{\Rbar \leftarrow P}$}, since the algorithm commutes
  with base change (fixing the $\{X_i(t_i), Y_i(t_i)\}$), and since
  \mbox{$\Ties_{\Rbar \leftarrow P}$} is unique as a subspace of
  \mbox{$\Tisec_{\Rbar \leftarrow P}$}.
In this way, we obtain an effective way to compute \mbox{$\Ties_{\Rbar
    \leftarrow  P}$} and, hence, to compute the semiuniversal equisingular
deformation of \mbox{$\Rbar \leftarrow P$} (see Theorem
\ref{thm:Defes}). \hfill\qedsymbol 
\end{rem}

\noindent
Proposition \ref{prop:ES=HN}, together with
the relation between (equisingular) deformations of the parametrization and
(equisingular) deformations of the equation discussed in Theorem
\ref{thm:fitting} and Theorem \ref{thm:esdef param, eqn} leads to the following
algorithm for computing the 
{\it equisingularity stratum\/} in the
base space \mbox{$A=K[[\bs]]/I$} of a deformation with trivial section of a reduced 
plane curve singularity (given by \mbox{$F\in K[[\bs,x,y]]$},
\mbox{$\bs=(s_1,\dots,s_N)$}). That is, the algorithm computes an ideal
\mbox{$\ES(F)\subset A$} such that the induced deformation over $A/\ES(F)$ is
equisingular along the trivial section and $\ES(F)$ is minimal in the
sense that, for each ideal \mbox{$J\subset A$} such that the induced 
deformation over $A/J$ is equisingular along the trivial section, we
have \mbox{$\ES(F)\subset J$}. 

\begin{alg}[Equisingularity stratum] 
\label{alg:2}\hfill

\begin{entry}
\item[Input] \mbox{$F\in K[[\bs,x,y]]$}, \mbox{$\bs=(s_1,\dots,s_N)$}, defining a
  deformation over the local $K$-algebra \mbox{$A=K[[s_1,\dots,s_N]]/I$} of the
  reduced plane curve singularity $C$ with equation \mbox{$f=F\;\mod \fm_A$}.
\item[Assume] Either \mbox{$\charac(K)=0$} or \mbox{$\charac(K)>\ord(f)$}.
\item[Output] A set of generators for \mbox{$\ES(F)\subset A$}.  
\end{entry}

\medskip\noindent
{\it Step 1.}\, {\it (Initialization)\/}
\smallskip
\begin{itemize}
\itemsep3pt
\item Compute the system $\kh^{(1)}_K, \dots,\kh^{(r)}_K$ of
  Hamburger-Noether expansions for \mbox{$f\in
    K[[x,y]]$}.\footnote{This may be done by applying
    the algorithm of Rybowicz \cite{Rybowicz} (extending the
    algorithm in \cite{Cam} to the reducible case). An implementation
    of this algorithm is provided by the {\sc Singular} library
    \texttt{hnoether.lib} written by M.\ Lamm.} In particular,
  determine the number $r$ of 
    branches of $C$.  
\item Set \mbox{$\kg\!\!\::=\!\!\:\emptyset$},\:
  \mbox{$n\!\!\::=\!\!\:\ord(f)$}.
\item For each \mbox{$i=1,\dots, r$}, set 
\mbox{$e[i]\!\!\::=\!\!\:F[i]\!\!\::=\!\!\:ok[i]\!\!\::=\!\!\:0$}.
\end{itemize}

\medskip\noindent
{\it Step 2.}\, {\it (Check equimultiplicity)\/}
\smallskip
\begin{itemize}
\itemsep3pt
\item If \mbox{$n=1$} then {\sc Return}$(\kg)$\,. 
\item Let \mbox{$F=\sum_{(\alpha,\beta)}
    a_{\alpha\beta} x^\alpha y^\beta$} then set
$$ \kg:= \kg\cup \{ a_{\alpha\beta}\mid \alpha+\beta<n\} \,,
\quad F:=F-\!\sum_{\alpha+\beta<n}\!  a_{\alpha\beta} x^\alpha
y^\beta\,.$$ 
\item Let the $n$-jet of \mbox{$f$} decompose as 
$$f  \equiv \overline{c}\cdot
x^{n_1}\!\cdot \prod_{\nu=2}^{\rho}
(y-\overline{a}_{\nu}x)^{n_{\nu}}\:\mod \langle
x,y\rangle^{n+1}\,,\quad \overline{a}_{\nu}\neq \overline{a}_{\nu'} \text{ for 
  } \nu\neq \nu'\,,$$
where the factor $x^{n_1}$ corresponds to $r_1$ branches of $C$, say
$C_1,\dots,C_{r_1}$,
while each factor \mbox{$(y-\overline{a}_{\nu}x)^{n_{\nu}}$},
\mbox{$\nu=2,\dots, \rho$}, 
corresponds to \mbox{$r_{\nu}-r_{\nu-1}$} bran\-ches, say 
$C_{r_{\nu-1}+1},\dots,C_{r_{\nu}}$
(this information can easily be read
from the Ham\-bur\-ger-Noether expansions $\kh^{(1)}_K, \dots,\kh^{(r)}_K$). 
Then we introduce new
variables \mbox{$b_{1},\dots,b_{\rho}$} and impose the
following 
condition on the $n$-jet (in $x,y$) of $F$:
\begin{equation}\label{eq:Step 5}
 \sum_{\alpha+\beta=n}  a_{\alpha\beta} x^\alpha
y^\beta \stackrel{!}{=} c\cdot
(x-b_1y)^{n_1}\!
\cdot \prod_{\nu=2}^{\rho}
\bigl(y-(b_{\nu}\!\!\:+\!\!\;\overline{a}_{\nu})x\bigr)^{n_{\nu}}
\end{equation}
with \mbox{$c\in A^\ast$}, \mbox{$c\equiv \overline{c} \ \mod \fm_A$}.
Set \mbox{$\rho_0:=\rho$}, and add the conditions
obtained by comparing the (\mbox{$n+1$}) coefficients of
\mbox{$x^{\alpha}y^{\beta}$}, \mbox{$\alpha+\beta=n$}, on both sides
of the equation to $\kg$. Note that $\kg$ is now a subset of
\mbox{$A[[b_1,\dots,b_{\rho_0}]]$}. 
\end{itemize}

\medskip\noindent
{\it Step 3.}\, {\it (1st blowing up)\/} 

\smallskip\noindent
 If \mbox{$r_1\!>0$} then set
  \mbox{$F[1]:=F(yx\!\!\:+\!\!\:b_1x,x)/x^{n}$}, \mbox{$n[1]:=n_1$}. Moreover,
  set  
$$F[r_{\nu}\!+1]:= \frac{F(x,yx\!\!\:+\!\!\:b_{\nu}x\!\!\:+\!\!\:
\overline{a}_{\nu}x)}{x^{n}}\,,\quad n[r_{\nu}\!+1]:= n_{\nu}\,,$$
\mbox{$\nu=1,\dots,\rho_0-1$}.

\medskip\noindent
{\it Step 4.}\, {\it (Check equimultiplicity after successive blowing
    up)\/} 

\smallskip\noindent
While \mbox{$S:=\{ i \:|\:F[i]\neq 0 \text{ and }
    ok[i]\neq 1\}\neq \emptyset$}, choose any \mbox{$i_0\in S$} and  
   do the following:

\smallskip
\begin{itemize}
\itemsep3pt
\item Set \mbox{$f[i_0]:=F[i_0]\:\mod \fm_A$}, and \mbox{$n:=\ord
    f[i_0]$}.  

\item If \mbox{$e[i_0]>1$} then the $n$-jet of $f[i_0]$
  necessarily equals
$y^{n}$, and we impose the following condition on the $n$-jet of
$F[i_0]$: 
\begin{equation}
  \label{eq:Step7b}
 F[i_0] \stackrel{!}{\equiv} c\cdot y^{n} \:\mod \langle
x,y\rangle^{n+1}\,.  
\end{equation}
Set \mbox{$e[i_0]:=e[i_0]-1$}, and add the conditions
obtained by comparing the co\-effi\-cients of
\mbox{$x^{\alpha}y^{\beta}$}, \mbox{$\alpha+\beta=n$}, on both sides
of the equation \eqref{eq:Step7b} to $\kg$. Finally, set
\mbox{$n[i_0]:=n$}, reduce $F[i_0]$ by the linear elements of $\kg$, and set
$$F[i_0]:=\frac{F[i_0](x,yx)}{x^{n}}\,.$$

\item Otherwise, redefine $\rho,n_{\nu},r_{\nu},\overline{a}_{\nu}$ such that 
$$f[i_0]  \equiv \overline{c}\cdot
x^{n_1}\!\cdot \prod_{\nu=2}^{\rho}
(y-\overline{a}_{\nu}x)^{n_{\nu}}\:\mod \langle
x,y\rangle^{n+1}\,,\quad \overline{a}_{\nu}\neq \overline{a}_{\nu'} \text{ for   
  } \nu\neq \nu'\,,$$
where the factor $x^{n_1}$ corresponds to $r_1$ branches, say 
$C_{i_0},\dots,C_{i_0+r_1-1}$, while each factor
\mbox{$(y-\overline{a}_{\nu}x)^{n_{\nu}}$},  \mbox{$\nu=2,\dots, \rho$},
corresponds to \mbox{$r_{\nu}-r_{\nu-1}$} bran\-ches, say
$C_{i_0+r_{\nu-1}},\dots,C_{i_0+r_{\nu}-1}$ (again, this information can easily
be read from $\kh^{(1)}_K, \dots,\kh^{(r)}_K$). We
introduce variables \mbox{$b_{\rho_0+1},\dots,b_{\rho_0+\rho-1}$} and impose the
following condition on the $n$-jet of $F[i_0]$:
\begin{equation}
  \label{eq:Step7c}
 F[i_0] \stackrel{!}{\equiv} c\cdot x^{n_1}\!
\cdot \prod_{\nu=2}^{\rho}
\bigl(y-(b_{\rho_0+\nu-1}\!\!\:+\!\!\;\overline{a}_{\nu})x\bigr)^{n_{\nu}}\:\mod
\langle 
x,y\rangle^{n+1}
\end{equation}
with \mbox{$c\in A^{\ast}$}, \mbox{$c\equiv \overline{c} \ \mod \fm_A$}.
Set \mbox{$\rho_0:=\rho_0+\rho-1$}, and add the conditions
obtained by comparing the coefficients of
\mbox{$x^{\alpha}y^{\beta}$}, \mbox{$\alpha+\beta=n$}, on both sides
of \eqref{eq:Step7c} to $\kg$. Reduce $F[i_0]$ by the linear elements
of $\kg$.  

\item {\it(Blowing up)} 

\smallskip\noindent
For \mbox{$\nu=\rho-1,\dots, 2$}, set 
$$F[i_0\!\!\:+r_{\nu}]:= \frac{F[i_0](x,yx\!\!\:+\!\!\:b_{\rho_0+\nu-1}x\!\!\:+\!\!\:
\overline{a}_{\nu}x)}{x^{n}}\,, \quad n[i_0\!\!\:+r_{\nu}]:=n_{\nu}\,.$$
Moreover, if \mbox{$r_1> 0$} then set \mbox{$F[i_0]:=F[i_0](yx,x)/x^{n}$},
$$e[i_0]:=\left\lceil \frac{n[i_0]\!\!\:- \!\!\:n_2\!\!\:-
  \!\!\:\ldots\!\!\:- \!\!\:n_{\rho}}{n_1}\right\rceil-1\,,$$ 
and \mbox{$n[i_0]:=n_1$}. 
\item If \mbox{$\ord F[i_0]\leq 1$} and \mbox{$e[i_0]\leq 1$} then
  \mbox{$ok(i_0):=1$}. 
\end{itemize}

\medskip\noindent
{\it Step 5.}\, {\it (Eliminate auxiliary variables)}

\smallskip\noindent
\begin{itemize}
\itemsep3pt
\item Set \mbox{$B:=\{1,\dots,\rho_0\}$}.
\item For each \mbox{$k\in B$} check whether in $\kg$ there is an element of type
  \mbox{$ub_{k}-a$} with \mbox{$u\in 
    A^\ast$}, \mbox{$a\in A[[\bb\setminus \{b_k\}]]$}. If yes, then
  replace \mbox{$b_{k}$} by  
  \mbox{$a/u\in A[[\bb\setminus \{b_k\}]]$} in all terms of elements of $\kg$, and set
  \mbox{$B:=B\setminus \{k\}$}.\,\footnote{This step applies, in particular, to all
    those $b_{k}$ which were introduced in an equation (\ref{eq:Step
      5}), resp.\ (\ref{eq:Step7c}), with $f[i_0]$
    being unitangential (see Remark \ref{rmk:poly data}).}
\item ({\it Hensel lifting}) For the remaining \mbox{$k\in B$} do the following: if
  $b_k$ appears only in one element of $\kg$, remove this
  element from $\kg$. Otherwise, compute the
  unique Hensel lifting of the factorization of
  \mbox{$(F[i_0]\:\mod \fm_A)\big|_{x=1}$} in the defining equation
  (\ref{eq:Step 5}), resp.\ (\ref{eq:Step7c}), for $b_k$:
$$  F[i_0](1,y) \equiv c\cdot  \prod_{\nu=2}^{s} g_{\nu} \ \mod
\langle x,y\rangle^{n+1}\,, \quad \
g_{\nu} \equiv (y-\overline{a}_{\nu})^{n_{\nu}}\:\mod \fm_A\,,
$$
where \mbox{$c\in A^\ast$}, and
\mbox{$g_{\nu}=y^{n_{\nu}}\!+c_{\nu}y^{n_{\nu}-1}\!+(\text{lower terms in $y$})\in
  A[y]$}. If the auxiliary variable $b_{k}$ was introduced in the factor
with constant term $\overline{a}_{\nu}^{n_{\nu}}$, then replace $b_{k}$ by
\mbox{$-(c_{\nu}/n_{\nu})-\overline{a}_{\nu}\in A$} in  all terms
of elements of $\kg$.\footnote{Note that, if the
  Hensel lifting for the factorization of \mbox{$f[i_0]\big|_{x=1}$} in the
  defining equation (\ref{eq:Step 5}) has to be computed, and if there
  is one factor of $f[i_0]$ with tangent $x$ and one with tangent 
$y$, apply a coordinate change of type
\mbox{$(x,y)\mapsto (x\!\!\:+\!\!\:\eta y,y)$}, \mbox{$\eta\in K$}, first.} 
\end{itemize}

\medskip\noindent
{\it Step 6.}\, {\sc Return}($\kg$). 
\end{alg}

\noindent
The proof of correctness for this algorithm is based on results of
\cite{Cam1} and the following two easy lemmas (see the end of this
section for proofs):

\begin{lem}[Uniqueness of Hensel lifting]\label{lem:unique Hensel}
Let \mbox{$A=\K[[t_1,\dots,t_r]]/I$} be a complete local $\K$-algebra, 
and let \mbox{$F\in A[y]$} be a monic polynomial satisfying 
$$ F\equiv (y+\overline{a}_1)^{m_1} \!\cdot \ldots \cdot
(y+\overline{a}_s)^{m_s} \:\modulo\, \fm_A\,,\quad \overline{a}_i\neq
\overline{a}_{i'}\in K \text{ for } i\neq i'\,. $$
Then there exists a unique Hensel lifting of the factorization,
$$F=g_1\cdot \ldots \cdot g_s \,,\quad g_i\in A[y] \text{ monic}\,,
\ g_i\equiv (y+\overline{a}_i)^{m_i} \:\modulo\, \fm_A\,.$$
\end{lem}

\begin{lem}\label{lem:unique section}
  Let $A$ be a local $K$-algebra,
  and suppose that the characteristic of $\K$ does not divide the
  positive integer $m$. Then, for any 
  \mbox{$a,b\in A$}, the following are equivalent: 
  \begin{enumerate}
  \itemsep0pt 
  \item[(1)] \mbox{$(y+a)^m= (y+b)^m\in A[[x,y]]$}\,,
 \item[(2)] $a=b$\,.
  \end{enumerate}
\end{lem}

\noindent
As mentioned before, the algorithm is based on the relation between equisingular
deformations of the equation (along the trivial section) and
Hamburger-Noether deformations. It is not difficult to see that the
terms \mbox{$(b_{\rho_0+\nu-1}\!\!\:+\!\!\;\overline{a}_{\nu})$} on  
the right-hand side of (\ref{eq:Step 5}), respectively
(\ref{eq:Step7c}), correspond precisely to the `free' coefficients
\mbox{$a_{j,k}^{(i)}$} of the Hamburger-Noether expansions
$\kh_A^{(i)}$, respecting the condition (HN2). The condition that the 
first nonzero coefficient in each row (except in the first one)
has to be a unit is reflected in the algorithm by introducing
$e[i_0]$. 
On the other hand, the left-hand side of (\ref{eq:Step 5}),
resp.\ (\ref{eq:Step7c}), is the deformation of $f$ obtained after 
performing the respective blowing-ups (with indeterminates $b_{\nu}$).
The proof of \cite[Thm.\ 1.3]{Cam1} shows that 
\mbox{$F$} defines an equisingular deformation of $R=P/\langle
f\rangle$ over 
$A/J$ along the trivial section $\sigma$ iff 
it defines an equimultiple deformation along $\sigma$ (Step 
2) and there exist \mbox{$b_{k}=b_{k}(\bs)\in A$},
\mbox{$k=1,\dots,\rho_0$}, such that the conditions (\ref{eq:Step 5}), 
(\ref{eq:Step7b}) and (\ref{eq:Step7c}) are satisfied modulo $J$.

Lemma \ref{lem:unique Hensel} implies that the factor
\mbox{$(y-(b_{\rho_0+\nu-1}\!\!\:+\!\!\;\overline{a}_{\nu})x)^{n_{\nu}}\in
  A[x,y]$} on the right-hand side of (\ref{eq:Step 5}), resp.\
(\ref{eq:Step7c}), is uniquely determined (as a factor of the Hensel
lifting of the factorization of \mbox{$f[i_0]=F[i_0]\:\mod
  \fm_A$}). Lemma \ref{lem:unique 
  section}, together with our assumption on the characteristic of $K$,
gives that $b_{\rho_0+\nu-1}$ is uniquely determined (as described in
Step 5 of the algorithm). Note that the integer $n_\nu$ appearing in
the Hensel lifting step of the algorithm is the sum of multiplicities of
the strict transforms of some branches of $C$, hence \mbox{$n_\nu\leq
\ord(f)$} and our assumption implies that $n_\nu$ is not divisible by the
characteristic of $K$. 
\phantom{a}\hfill\qedsymbol

\medskip
\begin{rem}[Working with polynomial data]\label{rmk:poly data}
In practice, we want (and can) apply Algorithm \ref{alg:2} only to the case
where the curve $C$ and its deformation are given by polynomials. 
Thus, let \mbox{$A=K[[\bs]]/I_0K[[\bs]]$} for some ideal \mbox{$I_0\subset 
  K[\bs]$}, and let \mbox{$F\in K[\bs,x,y]$}.  
Applying Algorithm \ref{alg:2} to $F$ does not necessarily lead 
  to polynomial (representatives of) generators for \mbox{$\ES(F)\subset A$}. 
This is caused by the Hensel lifting in Step 5. However,
under certain circumstances the Hensel lifting may be avoided, replacing Step
  5 by a Gr\"obner basis computation\footnote{A {\sc Singular}
implementation of the resulting algorithm is accessible via the command
\texttt{esStratum} provided by the library
\texttt{equising.lib} \cite{LoM}.}:

\medskip\noindent
{\it Step 5'.}\, {\it (Eliminate $\bb=(b_1,\dots,b_{\rho_0})$)}

\begin{itemize}
\item[]Let $J$ be the ideal of \mbox{$(K[\bs]_{\langle
      \bs\rangle}/I_0K[\bs]_{\langle \bs\rangle})[\bb]$} generated 
  by $\kg$. Compute a set of polynomial generators $\kg'$ for the elimination
  ideal $$J\cap (K[\bs]_{\langle \bs\rangle}/I_0K[\bs]_{\langle
      \bs\rangle})\,.$$
  This can be done by computing a
    Gr\"obner basis for $J$ with respect to a product ordering
    $(>_{\bb},>_{\bs})$ on $K[\bb,\bs]$, where $>_{\bb}$ is global and
    $>_{\bs}$ is local. Set \mbox{$\kg:=\kg'$}.
\end{itemize}

\noindent
Let, for instance, \mbox{$f=F(x,y,0)$} define an {\em irreducible\/} plane
curve singularity. Then all appearing polynomials 
\mbox{$F[i_0]\:\mod\fm_A$} are unitangential. Hence,
  (\ref{eq:Step7c}) reads either \mbox{$F[i_0]\equiv c\cdot x^{n}\:\mod
    \langle x,y\rangle^{n+1}$}, or
$$\begin{array}{rcl}
 F[i_0] &\equiv& c \cdot 
(y-(b_k\!\!\:+\!\!\;\overline{a})x)^{n} \nonumber \\ 
&\equiv& c\cdot (y^{n}\!-n (b_k\!\!\:+\!\!\;\overline{a})x
y^{n-1}\!+ x^2\!\cdot h(x,y))
\:\mod \langle x,y\rangle^{n+1}\,.   \label{eq:Step7c special}
\end{array}
$$
If the $n$-jet of $F[i_0]$ is \mbox{$\sum_{\alpha+\beta=n}
  a_{\alpha,\beta} x^\alpha y^\beta$} then the latter gives the
equations
\begin{equation}
  \label{eq:irred}
 c=a_{0,n}\in K[\bs]\setminus \langle \bs\rangle \,,\qquad nc\cdot
 b_k= -a_{1,n-1}-nc\overline{a} \in A\,.     
\end{equation}
In particular, the substitution of $b_k$ by 
\mbox{$-a_{1,n-1}/nc -\overline{a}$} in the elements of $\kg$ is
also performed by the Gr\"obner basis algorithm (multiplying the
resulting elements by appropriate units of the local ring
$K[\bs]_{\langle \bs\rangle}$). 

Similarly, if we consider a deformation over an Artinian base space,
say \mbox{$A=K[\bs]/\langle \bs\rangle^N$}, then we may again replace
Step 5 in the algorithm by the above Step 5'. In this case, we
additionally have to add to $\kg$ all monomials in $\bs,\bb$ of degree
$N$. 

In particular, this allows us to compute a set of generators for 
Wahl's equi\-sin\-gu\-larity ideal \cite{Wah} working with polynomial data only:

\begin{alg}[Equisingularity ideal] 
\label{alg:3}\hfill

\begin{entry}
\item[Input] \mbox{$f\in K[x,y]$}, defining a reduced plane curve
  singularity $C$.
\item[Assume] Either \mbox{$\charac(K)=0$} or \mbox{$\charac(K)>\ord(f)$}.
\item[Output] A set of generators for the equisingularity ideal
$$ \IES(f) := \left\{ g\in K[[x,y]] \left| \begin{array}{c}
      f+\varepsilon g \text{ defines an equisingular}\\
  \text{deformation of $C$ over }\Keps
\end{array}
\right.\right\}.$$  
\end{entry}

\medskip\noindent
{\it Step 1.}\, {\it (Initialization)\/}

\smallskip
\begin{itemize}
\itemsep3pt
\item Compute a (monomial) $K$-basis
  \mbox{$\{g_1,\dots,g_N\}\subset K[x,y]$} for 
  the $K$-al\-ge\-bra \mbox{$ \langle x,y\rangle \cdot \K[x,y]/(\langle
  f\rangle + \langle x,y\rangle \cdot \langle \frac{\partial
    f}{\partial x},\frac{\partial f}{\partial x}\rangle)$} (see \cite{GrP}).
\item Compute the system $\kh^{(1)}_K, \dots,\kh^{(r)}_K$ of
  Hamburger-Noether expansions for \mbox{$f\in K[[x,y]]$}. In
  particular, read the number $r$ of branches of $C$ 
  and the number $\rho_0$ of free
  infinitely near points of $C$ corresponding to non-nodal
  singularities of the reduced total transform 
  of $C$.
\item Introduce new variables $b_1,\dots,b_{\rho_0}$ and set
$$\kg\!\!\::=\!\!\:\{s_js_{j'},b_kb_{k'},s_jb_k \mid
    1\leq j,j'\leq N, 1\leq k,k'\leq \rho_0\}\subset K[\bs,\bb]\,,$$
  \mbox{$n\!\!\::=\!\!\:\ord(f)$}.
\item For each \mbox{$i=1,\dots, r$}, set 
\mbox{$e[i]\!\!\::=\!\!\:F[i]\!\!\::=\!\!\:ok[i]\!\!\::=\!\!\:0$}.
\end{itemize}

\item[{\it Step 2--4.}] As in Algorithm \ref{alg:2}, applied to
\mbox{$F = f + \sum_{k=1}^N s_k g_k \in K[\bs,x,y]$} and the ring
\mbox{$A=K[[\bs]]$}, \mbox{$\bs=(s_1,\dots,s_{N})$}. (Instead of
introducing new variables $b_k$, reuse the
variables $b_1,\dots,b_{\rho_0}$ introduced in Step 1).

\medskip\noindent
{\it Step 5'.}\, As above.

\medskip\noindent
{\it Step 6.}\, Compute a reduced normal form for $F$ w.r.t.\
  $\langle \kg'\rangle$ and set
$$ \kf = \bigl\{ F|_{\bs=e_i}\!-f\:\big|\;i=1,\dots,N\bigr\} \cup \bigl\{
f,\,\tfrac{\partial f}{\partial x},\, \tfrac{\partial f}{\partial
  y}\bigr\}\,.$$

\medskip\noindent
{\it Step 7.}\, {\sc Return}($\kf$).
\end{alg}

\noindent
A {\sc Singular} implementation of this algorithm is accessible via the
\texttt{esIdeal} command provided by \texttt{equising.lib} \cite{LoM}. 
  
Finally, also for reducible plane curve singularities, we may replace the
Hensel lifting step by Step 5'. Then the algorithm computes 
defining equations for the  equisingularity
($\mu$-constant) stratum as an algebraic subset
of \mbox{$V(I)\subset\Spec K[[\bs]]$} (but not necessarily with the
correct scheme-theoretic structure imposed by deformation
theory). Indeed, the computation in 
Step 5' yields equations for the image of $V(\kg)$ under the projection 
$$ \pi: \A^{\rho_0}\! \times (V(I_0), 0)^{\text{alg}}\to
  (V(I_0), 0)^{\text{alg}}\,,$$ where $(V(I_0), 0)^{\text{alg}}$ denotes
the germ of $V(I_0)$ at the origin with 
respect to the Zariski topology (see \cite{GrP}). Now, $V(\kg)$ intersects the
Zariski closure of the fibre \mbox{$\pi^{-1}(0)$} in 
\mbox{$\P^{\rho_0}\!\times \{0\}$} only
at the origin $\bo$ and at finitely many points \mbox{$\overline{\bb}$} which
correspond to a permutation of the factors in (\ref{eq:Step7b}),
resp.\ in (\ref{eq:Step7c}) (that is,
$b_{\rho_0+\nu-1}$ is replaced by 
\mbox{$b_{\rho_0+\nu'-1}\!\!\:+\!\!\;\overline{a}_{\nu'}
\!\!\:-\!\!\;\overline{a}_{\nu}$}) . The uniqueness of the Hensel lifting
implies that the image of the analytic germ of $V(\kg)$ at $\overline{\bb}$ under
$\pi$ coincides with the image of the analytic germ of $V(\kg)$ at $\bo$.   
Thus, the analytic germ of the image computed by eliminating $\bb$
coincides with the image under $\pi$ of the analytic germ of $V(\kg)$ at the origin.
\phantom{a}\hfill \qedsymbol
\end{rem}

\smallskip\noindent
We close this section by giving the proofs of Lemmas \ref{lem:unique
  Hensel}, \ref{lem:unique section}.

\smallskip\noindent
{\it Proof of Lemma \ref{lem:unique Hensel}.}
The existence of the Hensel lifting follows since $K[[\bs]]$ is
Henselian (see, e.g., \cite[\S\:I.5, Satz 6]{GrR}). It remains to prove
the uniqueness. Consider  
$$ G(y):=F(y-\overline{a}_s)=
\underbrace{g_1(y-\overline{a}_s)\cdot \ldots \cdot
  g_{s-1}(y-\overline{a}_s)}_{\textstyle =:u} \cdot
\underbrace{g_s(y-\overline{a}_s)}_{\textstyle =:h}\,, $$  
where \mbox{$u(0)\equiv \prod_{i=1}^{s-1}
  (\overline{a}_i-\overline{a}_s)\not\equiv 0$}. 
Hence, \mbox{$u\in A[[y]]^\ast$}, while $h$ is a Weierstra\ss\
polynomial in $A[y]$ (of degree $m_s$). 
Assuming that there exist two
such decompositions \mbox{$G=uh=u_1h_1$}, we would have
\mbox{$ 0=G\cdot (u^{-1}\!-\!\!\;u_1^{-1})+r\!\!\:-\!\!\:r_1$}, where
\mbox{$r\!\!\:-\!\!\:r_1\in A[y]$} has degree at most \mbox{$m_s\!-1$}. But 
$$G\,\equiv\,
  \overline{c}y^{m_s}\!+\text{\:(higher terms in $y$)} \:\mod \fm_A\,,\qquad
\overline{c}\in K\setminus\{0\}\,,$$ 
whence $G$ contains a term \mbox{$cy^{m_s}$}
(\mbox{$c\in A^\ast$}). 
Setting \mbox{$u^{-1}\!-\!\!\;u_1^{-1} =: \sum_{\alpha} c_{\alpha}
  y^\alpha$}, and 
choosing \mbox{$m\geq 0$} minimally such that \mbox{$C_m:= \left\{
    \alpha \geq 0 \:\big|\:  
    c_{\alpha}\in \fm_A^m\setminus \fm_A^{m+1} \right\} \neq
  \emptyset$}, it is obvious that the product \mbox{$cy^{m_s}\cdot
  c_{\alpha}y^\alpha\neq 0$} (with \mbox{$\alpha\in C_m$} minimal)
would have degree at least $m_s$, and could not be cancelled by any other
term of \mbox{$G\cdot (u^{-1}\!-\!\!\;u_1^{-1})$}. It follows that \mbox{$u=u_1$},
hence the uniqueness. 
\hfill\qedsymbol

\medskip\noindent
{\it Proof of Lemma \ref{lem:unique section}.}
In characteristic zero, the equivalence is obvious. Thus, let
\mbox{$\charac(\K)=p>0$} and write 
\mbox{$m=p^j\!\cdot\overline{m}$}, with $j$ a non-negative integer,
such that $p$ does not divide $\overline{m}$. Then the equality
\mbox{$(y+a)^m= (y+b)^m$} implies 
$$ \renewcommand{\arraycolsep}{3pt} 
\begin{array}{rcl}
0& = &(y+a)^{p^j\overline{m}}-(y+b)^{p^j\overline{m}} =
\bigl(y^{p^j}\!\!+a^{p^j}\bigr)^{\overline{m}}\!- 
\bigl(y^{p^j}\!\!+b^{p^j}\bigr)^{\overline{m}} \\
&=&
\overline{m}\cdot \bigl(a^{p^j}\!\!-b^{p^j}\bigr)\cdot
y^{(\overline{m}-1)p^j} \!+ \:\text{lower terms in $y$}\\
& = & 
\overline{m}\cdot(a-b)^{p^j} \cdot 
y^{(\overline{m}-1)p^j} \!+ \:\text{lower terms in $y$}\,.
\end{array}
$$
Hence, if \mbox{$p^j\!=1$} (that is, if 
$p$ does not divide $m$) we get \mbox{$a=b$}. 
The proof shows that the equivalence in Lemma
\ref{lem:unique section} holds in arbitrary characteristic if the
$K$-algebra $A$ is reduced. 
\hfill \qedsymbol

\begin{rem}
(1) In concrete calculations, we have to distinguish carefully between
deformations which are equisingular along a given section and those
which are abstractly equisingular, that is, equisingular along some
section. The corresponding deformation functors are $\uDefes_{(C,0)}$
and $\uDefES_{(C,0)}$ as introduced in Definition \ref{def:DefES}.

\medskip\noindent
(2) Algorithm \ref{alg:2} computes the ideal $\ES(F)$ of the maximal
stratum in the parameter space such that the restriction of the family
defined by $F$ is {\em equisingular along the trivial section}. If a family
with non-trivial section $\sigma$ is given, then one has to trivialize
this section first and then to apply Algorithm \ref{alg:2} in order to
compute the stratum such that the 
family is {\em equisingular along $\sigma$}. For instance, the family given
by \mbox{$F=(x-s)^2+y^3$} is equisingular along the section
\mbox{$s\mapsto (s,0,s)$}, while Algorithm \ref{alg:2} computes
\mbox{$\ES(F)=\langle s\rangle$}, which means that \mbox{$\{0\}$} is
the maximal stratum of 
equisingularity along the trivial section \mbox{$s\mapsto (0,0,s)$}.

\medskip\noindent
(3) Let \mbox{$K=\C$} and let $F$ define the semiuniversal deformation
with (trivial) section 
of the reduced plane curve singularity $(C,0)$ given by \mbox{$f\in
  \C\{x,y\}$}, that is, \mbox{$F(x,y,\bs)=f(x,y)+\sum_{i=1}^N 
  s_ig_i(x,y)$}, where \mbox{$\{g_1,\dots,g_N\}\subset \C\{x,y\}$} represents
a $\C$-basis of \mbox{$\langle x,y\rangle \cdot \C\{x,y\}/(\langle
  f\rangle + \langle x,y\rangle \cdot \langle \frac{\partial
    f}{\partial x},\frac{\partial f}{\partial x}\rangle)$}. Then the ideal
$\ES(F)$ as computed by Algorithm \ref{alg:2}
defines the stratum of $\mu$-constancy along the trivial section of
the family defined by $F$. This stratum is isomorphic to the
$\mu$-constant stratum of the semiuniversal deformation of $(C,0)$ (without
section) given by  \mbox{$G(x,y,\bs)=f(x,y)+\sum_{i=1}^{\tau}
  s_ih_i(x,y)$}, where \mbox{$\{h_1,\dots,h_{\tau}\}\subset \C\{x,y\}$} represents
a $\C$-basis of the Tjurina algebra \mbox{$\C\{x,y\}/\langle
  f,\frac{\partial
    f}{\partial x},\frac{\partial f}{\partial x}\rangle$} (this
follows from Theorem \ref{thm:esdef param, eqn}). 

Note that the ideal $\ES(F)$ contains more information than just about the
$\mu$-constant stratum. It gives the semiuniversal equisingular family
such that every fibre has a singularity of Milnor number $\mu$ {\em at the
origin}.

\medskip\noindent
(4) The isomorphism between the $\mu$-constant strata in (3) is unique
on the tangent level and the corresponding tangent map
$$\Ties_f:=\Iesf(f)\big/(\langle
  f\rangle + \langle x,y\rangle \cdot \langle \tfrac{\partial
    f}{\partial x},\tfrac{\partial f}{\partial x}\rangle)
  \xrightarrow{\cong} \IES(f)\big/\langle 
  f,\tfrac{\partial f}{\partial x},\tfrac{\partial f}{\partial
    x}\rangle=:\TiES_f$$  is induced 
by the inclusion \mbox{$\langle x,y\rangle \hookrightarrow
  \C\{x,y\}$}. Here, 
$$ \Iesf(f) := \left\{ g\in K[[x,y]] \left| \begin{array}{c}
      f+\varepsilon g \text{ defines an equisingular deformation of}\\
  \text{$\{f=0\}$ over }\C[\varepsilon]\text{ along the
    trivial section }
\end{array} \right.\right\}\,,
$$
which can be computed along the lines of Algorithm \ref{alg:3}, replacing
the definition of $\kf$ in Step 6 by 
$$ \kf := \bigl\{ F|_{\bs=e_i}\!-f\:\big|\;i=1,\dots,N\bigr\} \cup \bigl\{
f,\,x\tfrac{\partial f}{\partial x},\, x\tfrac{\partial f}{\partial
  y}, y\tfrac{\partial f}{\partial x},\, y\tfrac{\partial f}{\partial
  y}\bigr\}\,.$$ 
The {\sc Singular} procedure \texttt{esIdeal} returns both, $\IES(f)$
and $\Iesf(f)$.
\end{rem}

\section{Examples}
\setcounter{thm}{0}
\setcounter{equation}{0}

\noindent 
In the first example, we compute defining equations for the
stratum of $\mu$-constancy along the trivial section for a deformation
of a reduced plane curve 
singularity (with two sin\-gu\-lar branches) over a smooth base. We
proceed along the lines of Algorithm \ref{alg:2}, slightly modifying
and anticipating Step 5 (resp.\ Step 5'):

\begin{exmp}
Let \mbox{$\charac(K)\neq 2$} and 
consider the deformation of the Newton degenerate plane curve singularity
\mbox{$C=\{(y^4\!+x^5)^2\!+x^{11}=0\}$} 
over \mbox{$A=K[[\bs]]$}, \mbox{$\bs=(s_1,\dots,s_{10})$}, given by 
$$\renewcommand{\arraycolsep}{3pt} 
\begin{array}{rcl}
F&:=&(y^4\!+x^5)^2\!+x^{11}+s_1x^3y^6\!+
s_2x^9y^3\!+s_3x^8y^3\!+s_4x^7y^3\!+
s_5x^{10}y^2\!\\
&&\phantom{(y^4\!+x^5)^2\!+x^{11}}+s_6x^9y^2\!+s_7x^8y^2\!+
s_8x^{10}y+s_9x^9y+s_{10}x^{10}\,.
\end{array}
$$
In the first step of the algorithm, we compute the system of
Hamburger-Noether expansions for $C$ (developing each final row up to a 
sufficiently high order as needed for computing the system of
multiplicity sequences):
$$
(\kh_A^{(1)}) \ \
\renewcommand{\arraystretch}{1.2} 
\renewcommand{\arraycolsep}{3pt} 
\begin{array}{rcl}
z_{-1}&=&z_0z_1\\
z_{0}&=&-z_1^4+z_1^6-\frac{3}{2}z_1^8+\ldots
\end{array}
\qquad
(\kh_A^{(2)}) \ \
\renewcommand{\arraystretch}{1.2} 
\renewcommand{\arraycolsep}{3pt} 
\begin{array}{rcl}
z_{-1}&=&z_0z_1\\
z_0&=&-z_1^4-z_1^6-\frac{3}{2}z_1^8+\ldots
\end{array}
$$
Since all deformation terms lie above (or on) the Newton
boundary, the equimultiplicity condition in Step 2 of the algorithm
does not lead to a new element of 
\mbox{$\kg$}. Further, we impose a
factorization \mbox{$y^8=c\cdot (y-b_1x)^8$}, which is only possible for
\mbox{$b_1=0$} (that is, \mbox{$\kg=\kg\cup\{b_1\}$}). We apply
the formal blowing-up (Step 3)
$$\renewcommand{\arraycolsep}{3pt} 
\begin{array}{rcl}
F[1]:=\dfrac{F(x,yx)}{x^8}&=&(y^4\!+x)^2\!+x^3+s_1xy^6\!+s_2x^4y^3\!+
s_3x^3y^3\!+s_4x^2y^3\\
&&\phantom{(y^4}\!+s_5x^4y^2+s_6x^3y^2\!+s_7x^2y^2\!+
s_8x^3y+s_9x^2y+s_{10}x^2
\end{array}
$$
and set \mbox{$n[1]:=8$}. We obtain
\mbox{$f[1]=(y^4\!+\!\!\;x)^2\!+x^3$} which has 
order \mbox{$n\!\!\;=\!\!\;2$}. Hence, in \eqref{eq:Step7c}, we impose
the condition  
\mbox{$F[1]\equiv cx^2\:\mod\langle x,y\rangle^3$}, which is obviously
satisfied for \mbox{$c=1+s_{10}\in A^\ast$}. We set 
$$\renewcommand{\arraycolsep}{3pt} 
\begin{array}{rcl}
F[1]:=\dfrac{F[1](yx,x)}{x^2}&=&(x^3\!+y)^2\!+xy^3+s_1x^5y+s_2x^5y^4\!+
s_3x^4y^3\!+s_4x^3y^2\!\\
&&\phantom{(y}\!\!\,s_5x^4y^4+s_6x^3y^3\!+s_7x^2y^2\!+
s_8x^2y^3+s_9xy^2\!+s_{10}y^2,
\end{array}
$$
\mbox{$e[1]:=\lceil 8/2 \rceil-1=3$} and \mbox{$n[1]:=2$}. Hence, in
the following two turns of the loop in Step 4, we impose the condition
\mbox{$F[1]\equiv cy^2\:\mod \langle x,y\rangle^3$} and perform then
the formal blowing-up \mbox{$F[1]:=F[1](x,yx)/x^2$}. Note that both turns do not
lead to new elements of $\kg$. After the second turn, we have 
$$
F[1]\equiv (x+y)^2\!+x^3y^3+s_1x^3y+
s_4x^3y^2\!+s_7x^2y^2\!+s_9xy^2\!+s_{10}y^2
$$ 
modulo \mbox{$\langle \kg\rangle +\langle
x,y\rangle^7$}. In the next turn, we impose the condition
$$
(1\!\!\:+\!\!\:s_{10})y^2+2xy+x^2 \stackrel{!}{=} c\cdot
(y-(b_2\!\!\:-\!\!\:1)x)^2 \,, 
$$
hence \mbox{$c=1\!\!\:+\!\!\:s_{10}$}, and we obtain the equations
$$  (1\!\!\:+\!\!\:s_{10})\cdot b_2-s_{10}=0\,,\qquad
(1\!\!\:+\!\!\:s_{10})\cdot (b_2-1)^2-1=0\,,$$ 
which imply \mbox{$b_2=s_{10}=0$}, that is, in Step 5 (resp.\ 5'),
$b_2$ and $s_{10}$ will be added to $\kg$ (we anticipate this here and set
\mbox{$\kg=\kg\cup \{b_2,s_{10}\}$}). 
We apply the formal blowing-up
$$\renewcommand{\arraycolsep}{3pt} 
\begin{array}{rcl}
F[1]&:=&\dfrac{F[1](x,yx-x)}{x^2}\\
&\,\equiv &y^2\!-x^4
+s_1(x^2y-x^2)+s_4(-2x^3y+x^3)+s_7(x^2y^2\!-2x^2y+x^2)\\
 &&\phantom{y^2\!-x^4}
+ s_9(xy^2\!-2xy+x)\:
\end{array}
$$
modulo \mbox{$\langle \kg\rangle +\langle x,y\rangle^5$}.
The imposed condition reads now
$$ y^2-2s_9xy+(s_7\!-\!\!\:s_1)x^2+s_9x \stackrel{!}{=} c\cdot
(y-b_3x)^2\,, 
$$
hence \mbox{$c=1$}, \mbox{$s_9=0$}, \mbox{$b_3=s_9$} and
\mbox{$s_7\!-\!\!\:s_1=b_3^2$}. That is, partly anticipating \mbox{Step 5}
or 5', we set \mbox{$\kg=\kg\cup
  \{s_9,b_3,s_7\!-\!\!\:s_1\}$}.  We apply the formal blowing-up 
$$
F[1]:=\frac{F[1](x,yx)}{x^2}\equiv y^2\!-x^2\!
-s_1xy+s_4x\ \, \mod \langle \kg\rangle +\langle x,y\rangle^3
$$
and impose the condition
$$ y^2\!-x^2\!-s_1xy+s_4x \stackrel{!}{=} c\cdot
(y\!\!\:-\!\!\:b_4x\!\!\:-\!\!\:x)(y\!\!\:-\!\!\:b_5x\!\!\:+\!\!\:x)\,.
$$
Hence, \mbox{$s_4=0$}, \mbox{$c=1$}, and we obtain the equations
(again partly anticipating Step 5 or 5'):
$$ b_4+b_5=s_1\,,\quad b_5^2-(2+s_1)b_5+s_1=0\,, $$
which are added to $\kg$. Now, we set 
$$ F[2]:=\frac{F[1](x,yx\!\!\:+\!\!\:b_5\!\!\:-\!\!\:x)}{x}\,,\qquad
F[1]:=\frac{F[1](x,yx\!\!\:+\!\!\:b_4\!\!\:+\!\!\:x)}{x}\,,$$ 
both being of order $1$, whence \mbox{$ok[1]=ok[2]=1$}, and we may
assume to enter Step 5 with
$$ \kg= \{ s_1\!-\!\!\:s_7,\, s_4,\, s_9,\, s_{10},\, b_1,\, b_2,\,
b_3,\,
b_4\!+\!\!\:b_5\!-\!\!\:s_1,\,
b_5^2\!\!\:-\!\!\:(2\!\!\:+\!\!\:s_1)b_5\!\!\:+\!\!\:s_1
\}\,.
$$
Since $b_4$ appears in exactly one of the elements of $\kg$, we simply
remove this element from $\kg$. Then $b_5$ appears in only one
element, too. So, we also remove this element and there is no need to
apply a Hensel lifting step, that is, to compute the power series
expansion of
\mbox{$b_5=\bigl(2+s_1-\sqrt{s_1^2+4}\bigr)\big/2$}. The same 
result is obtained by applying the elimination procedure of Step 5':
$$ ES(F) = \langle s_1\!-\!\!\:s_7,\, s_4,\, s_9,\, s_{10}\rangle
\subset \K[[\bs]]\,. $$
Since the deformation terms of $F$, together with the terms below the
Newton boundary, generate the Tjurina algebra \mbox{$\K[[x,y]]/\langle
  f,\frac{\partial f}{\partial x},\frac{\partial f}{\partial
    x}\rangle$},   
we can, in particular, read off the
equisingularity ideal of \mbox{$f=(y^4\!+x^5)^2\!+x^{11}$}:
$$ \IES(f)= \left\langle f,\,\tfrac{\partial f}{\partial x}\,,\,
  \tfrac{\partial f}{\partial y}\,,\, x^3y^6\!+x^8y^2\!,\, x^8y^3\!,\,
  x^9y^2\!,\, x^{10}y\right\rangle\,. 
$$

\vspace{-14pt}
\hfill\qedsymbol
\end{exmp}

\medskip\noindent
The second example shows the computation of $\IES(f)$ in the case of a 
Newton degenerate plane curve singularity with $8$ smooth
branches:

\begin{exmp}
Let \mbox{$f=(y^4\!-x^4)^2-x^{10}\in \K[x,y]$}.  We start with the
versal deformation with trivial section of $f$, given by \mbox{$F\in K[\bs,x,y]$},
\mbox{$\bs=(s_1,\dots,s_{50})$}, 
$$\begin{array}{rcl}
 F &=& (y^4\!-x^4)^2-x^{10}\!+s_1y^{11}+s_2xy^{10}\!+s_3y^{10}+s_4xy^9+
s_5y^9\!+s_6xy^8\\
&& \phantom{(x^4\!-y^4)^2-x^{10}\!}+s_7y^8\!+s_8x^3y^7+s_9x^2y^7\!+s_{10}xy^7
 \\
&& \phantom{(x^4\!-y^4)^2-x^{10}\!}
+s_{11}x^3y^6\!+s_{12}x^2y^6+s_{13}x^3y^5\!
+s_{14}x^6y^2\!+\ldots 
\end{array}
$$
(here, we displayed only the terms of degree at least $8$). The
system of Ham\-bur\-ger-Noether expansions for $f$ is 
$$
\renewcommand{\arraystretch}{2} 
\renewcommand{\arraycolsep}{3pt} 
\begin{array}{llll}
(\kh_A^{(1)}) \ \;
\renewcommand{\arraystretch}{1} 
\renewcommand{\arraycolsep}{3pt} 
\begin{array}{rcl}
z_{-1}&=&-z_0-\frac{1}{4}z_0^2+\ldots\\
\end{array}
 &\quad&
(\kh_A^{(2)}) \ \;
\renewcommand{\arraystretch}{1} 
\renewcommand{\arraycolsep}{3pt} 
\begin{array}{rcl}
z_{-1}&=&-z_0+\frac{1}{4}z_0^2+\ldots\\
\end{array}
\\
(\kh_A^{(3)}) \ \;
\renewcommand{\arraystretch}{1} 
\renewcommand{\arraycolsep}{3pt} 
\begin{array}{rcl}
z_{-1}&=&z_0+\frac{1}{4}z_0^2+\ldots\\
\end{array}
&\quad&
(\kh_A^{(4)}) \ \;
\renewcommand{\arraystretch}{1} 
\renewcommand{\arraycolsep}{3pt} 
\begin{array}{rcl}
z_{-1}&=&z_0-\frac{1}{4}z_0^2+\ldots\\
\end{array}
\\
(\kh_A^{(5)}) \ \;
\renewcommand{\arraystretch}{1} 
\renewcommand{\arraycolsep}{3pt} 
\begin{array}{rcl}
z_{-1}&=&iz_0-\frac{i}{4}z_0^2+\ldots\\
\end{array}
&\quad&
(\kh_A^{(6)}) \ \;
\renewcommand{\arraystretch}{1} 
\renewcommand{\arraycolsep}{3pt} 
\begin{array}{rcl}
z_{-1}&=&iz_0+\frac{i}{4}z_0^2+\ldots\\
\end{array}
\\
(\kh_A^{(7)}) \ \;
\renewcommand{\arraystretch}{1} 
\renewcommand{\arraycolsep}{3pt} 
\begin{array}{rcl}
z_{-1}&=&-iz_0+\frac{i}{4}z_0^2+\ldots\\
\end{array}
 &\quad&
(\kh_A^{(8)}) \ \;
\renewcommand{\arraystretch}{1} 
\renewcommand{\arraycolsep}{3pt} 
\begin{array}{rcl}
z_{-1}&=&-iz_0-\frac{i}{4}z_0^2+\ldots\\
\end{array}
\end{array}
$$
where \mbox{$i=\sqrt{-1}$}. From these expansions,
we read that there are $12$ free infinitely near points of
\mbox{$C=\{f=0\}$} corresponding to non-nodal singularities of the
reduced total transform of $C$. We initialize $\kg$ as
$$\kg\!\!\::=\!\!\:\{s_js_{j'},b_kb_{k'},s_jb_k \mid
    1\leq j,j'\leq 48, \,1\leq k,k'\leq 12\}\subset K[\bs,\bb]\,.$$
The equimultiplicity condition of Step 2 implies that
the 34 non-displayed terms of $F$ must be zero, that is, we set
\mbox{$\kg:= \kg\cup \{s_{15},\dots,s_{48}\}$}. 
In Step 5, we impose now a decomposition 
$$
\renewcommand{\arraycolsep}{3pt} 
\begin{array}{l}(y^4\!-x^4)^2+s_7y^8\!+s_{10}xy^7\!+s_{12}x^2y^6\!+s_{13}x^3y^5\!+
s_{14}x^6y^2\\
\qquad \stackrel{!}{=} c\cdot (y-b_1x-x)^2 \!\cdot(y-b_2x+x)^2 \!\cdot
(y-b_3x+ix)^2\!\cdot (y-b_4x-ix)^2 
\end{array}
$$
which modulo \mbox{$\langle \kg \rangle$} leads to 8 new linear relations:
$$\kg = \kg\cup \{ s_7, \,s_{10}, \,s_{13}, \,s_{12}\!+\!\!\;s_{14},
\, 8b_1\!-\!\!\:s_{14},\, 8b_2\!+\!\!\:s_{14},\, 8b_3\!-\!\!\:is_{14},
\,8b_4\!+\!\!\:is_{14}\}\,.$$
We set
$$ 
\renewcommand{\arraycolsep}{3pt} 
\begin{array}{rclcrcl}
F[1]&:=&\dfrac{F(x,yx+b_1x+x)}{x^2}\,,&\quad &
F[3]&:=&\dfrac{F(x,yx+b_2x-x)}{x^2} \,,\\[0.6em] 
F[5]&:=&\dfrac{F(x,yx+b_3x-ix)}{x^2}\,,&\quad &
F[7]&:=&\dfrac{F(x,yx+b_4x+ix)}{x^2}\,,
\end{array}
$$
all of them being of multiplicity \mbox{$2=n[1]=n[3]=n[5]=n[7]$} as
power series in $x,y$. Choosing, for  
instance, \mbox{$i_0=1$} (that is, considering $F[1]$), we impose
the new condition (modulo $\langle \kg\rangle$) 
$$
\renewcommand{\arraycolsep}{3pt} 
\begin{array}{l}
16y^2\!-x^2\!+(s_3+s_4+s_8)x^2\!+
(9s_5+8s_6+7s_9+6s_{11})xy\\
\qquad +\,4s_{14}y^2\!+(s_5+s_6+s_9+s_{11})x
\stackrel{!}{=}  c\cdot \bigl(y-b_5x+\tfrac{1}{4}x\bigr)
\cdot\bigl(y-b_6x-\tfrac{1}{4}x\bigr)\,, 
\end{array}
$$
which leads to the conditions
$$ 
\renewcommand{\arraycolsep}{3pt} 
\begin{array}{rcl}
J=J+\langle s_5+s_6+s_9+s_{11}&,&
32b_5\!\!\:-\!\!\:4s_3\!\!\:-\!\!\:4s_4\!\!\:-\!\!\:s_6\!\!\:-
\!\!\:4s_8\!\!\:-\!\!\:2s_9\!\!\:-\!\!\:3s_{11}\!\!\:-\!\!\:s_{14}\,,\\
&&32b_6\!\!\:+\!\!\:4s_3\!\!\:+\!\!\:4s_4\!\!\:-\!\!\:s_6\!\!\:+
\!\!\:4s_8\!\!\:-\!\!\:2s_9\!\!\:-\!\!\:3s_{11}\!\!\:+\!\!\:s_{14}
\,\rangle\,.
\end{array}
$$
Proceeding in the same way with the other possible choices
\mbox{$i_0=3,5,7$}, we obtain three more (linearly independent) conditions
for \mbox{$s_5, s_6, s_9, s_{11}$}, and conditions of type
\mbox{$b_k+L_k$}, $L_k$ some linear polynomial in $\bs$,\,
\mbox{$k=7,\dots,12$}. Since the eight polynomials
$F[1],\dots,F[8]$ obtained after the next formal blowing-ups 
are all of order $1$, we reach Step 5' and compute
$$ J \cap K[\bs]_{\langle \bs\rangle} = \langle s_5,\, s_6,\, s_7,\, s_9,\,s_{10},\,
s_{11}, \,s_{12}\!+\!\!\;s_{14},\,s_{13},
\,s_{15},\dots,s_{48}\rangle\,.$$
Hence, the base of the semiuniversal equisingular
deformation of $f$ has dimension $6$, and
$$ \IES(f)= \left\langle f,\,\tfrac{\partial f}{\partial x}\,,\,
  \tfrac{\partial f}{\partial y}\,,\, 
  x^6y^2 \!-x^2y^6\!,\,x^3y^7\!,\, xy^9\!,\,y^{10}
 \right\rangle\,.
$$

\vspace{-14pt}
\hfill\qedsymbol
\end{exmp}

\begin{rem} The correctness of the computed equations for the stratum
  of $\mu$-constancy (resp.\ equisingularity) can be checked by
choosing a random point $\bp$ satisfying the equations and computing
the system of Hamburger-Noether expansions for the evaluation
of $F$ at \mbox{$\bs=\bp$}. From the system of
Hamburger-Noether expansions, we can read a complete set of
numerical invariants of the equisingularity type (such as the
Puiseux pairs and the intersection numbers) which have to coincide
with the respective invariants of $f$. In characteristic $0$, it
suffices to compare the two Milnor numbers. We use {\sc
  Singular} to compute the $\mu$-constant stratum in our second example:

\begin{small}
\begin{verbatim}
   LIB "equising.lib";        //loads deform.lib, sing.lib, too
   ring R = 0, (x,y), ds;
   poly f = (y4-x4)^2 - x10;
   ideal J = f, maxideal(1)*jacob(f);
   ideal KbJ = kbase(std(J));
   int N = size(KbJ);
   ring Px = 0, (a(1..N),x,y), ls; 
   matrix A[N][1] = a(1..N);  
   poly F = imap(R,f)+(matrix(imap(R,KbJ))*A)[1,1]; 
   list M = esStratum(F);   //compute the stratum of equisingularity
                            //along the trivial section
   def ESSring = M[1]; setring ESSring;
   option(redSB);
   ES = std(ES);
   size(ES);                //number of equations for ES stratum
   //-> 42
\end{verbatim}
\end{small}

\noindent
Inspecting the elements of $\texttt{ES}$, we see that $40$ of the $48$
deformation parameters must vanish. Additionally, there are two
non-linear equations, showing that the equisingularity
($\mu$-constant) stratum is smooth (of dimension $6$) but not linear:

\begin{small}
\begin{verbatim}
   ES[1];
   //-> 8*A(1)+8*A(22)+A(1)^3
   ES[34];
   //-> 8*A(40)-A(1)^2+A(1)*A(22)
\end{verbatim}
\end{small}

\noindent
We reduce \texttt{F} by \texttt{ES} and evaluate the result at a
random point satisfying the above two non-linear conditions:
\begin{small}
\begin{verbatim}
   poly F = reduce(imap(Px,F),ES);  //A(1),A(22) both appear in F
   poly g = subst(F, A(22), -A(1)-(1/8)*A(1)^3); 
   for (int ii=1; ii<=44; ii++){ g = subst(g,A(ii),random(1,100)); }
   setring R;
   milnor(f);                       //Milnor number of f
   //-> 57
   milnor(imap(ESSring,g));         //Milnor number of g
   //-> 57
\end{verbatim}
\end{small}

\vspace{-14pt}
\hfill\qedsymbol
\end{rem}

\end{document}